\documentclass[preprint,12pt]{elsarticle}
\usepackage{amssymb}
\journal{Medical Image Analysis}
\usepackage{amsmath}
\usepackage{amssymb}
\usepackage{enumerate}
\usepackage{theorem}
\usepackage{graphicx,color}
\usepackage{pifont}
\usepackage{algorithm,algorithmic}
\newcommand{\eqdef}{\mbox{$\buildrel\triangle\over =$}}
\def\trans{\mbox{\tiny $\top$}}
\newcommand{\vect}[1]{{\boldsymbol #1}}    
\newcommand{\RR}{{\ensuremath{\mathbb R}}}

\newcommand{\NN}{{\ensuremath{\mathbb N}}}

\def\hermit{\mbox{\tiny\textsf{H}}}

\newcommand{\scal}[2]{\left\langle{{#1}\mid{#2}}\right\rangle} 
\newcommand{\KK}{\ensuremath{\mathbb K}}

\newcommand{\prox}{\ensuremath{\mathrm{prox}}}

\newcommand{\pinf}{\ensuremath{+\infty}}

\newtheorem{definition}{Definition}[section]
\newtheorem{theorem}{Theorem}[section]
\newtheorem{proposition}[theorem]{Proposition}
\theoremstyle{plain}{\theorembodyfont{\rmfamily}%
\newtheorem{example}[theorem]{Example}}

\newtheorem{assumption}[theorem]{Assumption}

\graphicspath{{./Figures/}}
\begin{document}

\begin{frontmatter}
\title{An Iterative Method for Parallel MRI SENSE-based Reconstruction in the Wavelet Domain}
 \author[label1,label2]{L. Cha\^ari}
 \author[label1]{J.-C. Pesquet}
 \author[label2]{Ph. Ciuciu}
 \author[label3]{and A. Benazza-Benyahia}
 \address[label1]{LIGM (UMR-CNRS 8049), Universit{{\'e}} Paris-Est, Champs-sur-Marne, 77454 Marne-la-Vall{\'e}e cedex, France. \\$\{$chaari,pesquet$\}$@univ-mlv.fr.}
 \address[label2]{CEA/DSV/$\mathrm{I}^2$BM/Neurospin, CEA Saclay, Bat. 145, Point Courrier 156, 91191 Gif-sur-Yvette cedex, France. \\ philippe.ciuciu@cea.fr}
 \address[label3]{Ecole Sup\'erieure des Communications de Tunis (SUP'COM-Tunis), URISA, Cit\'e Technologique des Communications, 2083, Tunisia.\\ benazza.amel@supcom.rnu.tn 
}
\address{\tnotetext[t1]}
\tnotetext[t1]{Part of this work was presented at the ISBI conference in 2008 \cite{chaari_08}.\\
This work was supported by R\'egion Ile de France, and the Agence Nationale de la Recherche under grant
ANR-05-MMSA-0014-01.}
\begin{abstract}
To reduce scanning time and/or improve spatial/temporal resolution in some MRI applications, parallel MRI (pMRI) acquisition techniques with multiple coils acquisition have emerged since the early 1990's as powerful 3D imaging methods that allow faster acquisition of reduced Field of View (FOV) images. In these techniques, the full FOV image has to be reconstructed from the resulting acquired undersampled $k$-space data. To this end, several reconstruction techniques have been proposed such as the widely-used
SENSE method. However, the reconstructed image generally presents artifacts when perturbations occur in both the measured data and the estimated coil sensitivity maps. In this paper, we aim at achieving good reconstructed image quality when using low magnetic field and high reduction factor. Under these severe experimental conditions, neither the SENSE method nor the Tikhonov regularization in the image domain give convincing results.
To this aim, we present a novel method for SENSE-based reconstruction which proceeds with regularization in the complex wavelet domain. To further enhance the reconstructed image quality, local convex constraints are added in the regularization process.
In vivo experiments carried out on Gradient-Echo (GRE) anatomical and Echo-Planar Imaging (EPI) functional MRI data at 1.5 Tesla indicate that the proposed algorithm provides reconstructed images with reduced artifacts for high reduction factor.
\end{abstract}


\end{frontmatter}
\section{Introduction}
Reducing global imaging time is of main interest in neuroimaging, in particular in the study of brain dynamics using functional MRI (fMRI). Actually, a short acquisition time lets one improve the spatial/temporal resolution of acquired fMRI data, which leads to a more efficient statistical analysis.
In addition, by reducing the global imaging time, some additional artifacts caused by the patient motion can be avoided.
For this reason, parallel imaging systems have been developed: multiple receiver surface coils with complementary sensitivity profiles located around the underlying object are employed to simultaneously collect in the frequency domain (i.e the $k$-space), data sampled at a rate $R$ times lower than the Nyquist sampling rate along at least one spatial direction, i.e. the phase encoding one. 
Therefore, the total acquisition time is $R$ times shorter than with conventional non parallel imaging. A reconstruction step is then performed to build a full Field of View (FOV) image by unfolding the undersampled measured ones. This reconstruction is a challenging task because of the low Signal to Noise Ratio (SNR) in parallel MRI (pMRI) caused by aliasing artifacts related to 
the undersampling rate, those caused by noise during the acquisition process and also the presence of errors in the estimation of coil sensitivity maps. The Simultaneous Acquisition of Spatial Harmonics (SMASH) \cite{Sodickson_D_97} was the first reconstruction method introduced by Sodickson and Manning in 1997, operating in the $k$-space domain. It uses a linear combination of pre-estimated coil sensitivity maps to generate the missing phase encoding steps. Some other $k$-space based reconstruction techniques have also been proposed like GRAPPA (Generalized Autocalibrating Partially Parallel Acquisitions) \cite{griswold_02}, a generalized implementation of variable-density AUTO-SMASH (VD-AUTO-SMASH) \cite{Heidemann_01}, which improves the AUTO-SMASH \cite{Jakob_98} method. The main difference between SMASH and these methods is that SMASH needs a separate coil sensitivity map estimation, while the other methods need the acquisition of few additional $k$-space lines to derive coil sensitivity maps without a reference scan. This is the reason that motivates the autocalibrating qualification of such methods. Indeed, GRAPPA may be preferred to non-autocalibrated methods when accurate coil sensitivity maps may be difficult to extract, or when reference scans are not possible to obtain at the right spatial resolution. For instance, in lung and abdomen imaging, the numerous inhomogeneous regions with a low spin density make inaccurate the estimation of the sensitivity information. However, all the reported methods may suffer from phase cancellation problems, low SNR during the acquisition process and limited reconstruction quality. 

In \cite{Pruessmann_K_99}, an alternative reconstruction method called SENSE (Sensitivity Encoding) has been introduced. It is a two-step procedure relying first on a reconstruction of reduced FOV images and second on a spatial unfolding technique, which amounts to a weighted least squares estimation. This technique requires a precise estimation of coil sensitivity maps using a reference scan (usually a 2D Gradient-Echo (GRE)).\\
To the best of our knowledge, in actual clinical daily routines, only GRAPPA and SENSE are offered by scanner manufacturers (or variants like mSENSE by Siemens). However, SENSE which is investigated in this work, 
remains the most frequently employed technique. In addition to neuroimaging, many other applications like cardiac imaging \cite{weiger_00} benefit from the enhanced image acquisition capabilities of this method. For a general overview of reconstruction methods in pMRI see \cite{Hoge_05}.

Actually, SENSE is often supposed to achieve an exact reconstruction in the case of non-noisy data and perfect coil sensitivity maps knowledge, which 
is also true for all above mentioned methods. However, in practice, noisy data and inaccuracies in the estimation of coil sensitivity maps make the reconstruction problem ill-conditioned. To overcome this limitation, some alternatives have been proposed like the optimization of the coil geometry \cite{Weiger_M_01} and the improvement of coil sensitivity profile estimation \cite{Blaimer_M_04}. But, when low magnetic field intensities (up to 1.5 Tesla) are used, high reduction factors are considered as unfeasible since the reconstructed images are affected by severe aliasing artifacts, and thus, a reduction factor value of two is often considered as the highest possible reduction factor at this field intensity. If one wants to further improve the spatial resolution to avoid filtering that blurs activation, to improve the temporal resolution in fMRI or to reduce the global imaging time in clinical MRI, it would be at the expense of the reconstructed image quality and/or potentially the reliability of statistical analysis in fMRI and brain activation detection.\\
As the image reconstruction problem is ill-posed, several works have dealt with regularization methods \cite{chaari_08,Ying_L_04,Rabrait07} to better estimate full FOV images, most of them operating in the image domain. For an introductory survey of linear inverse problems in application to pMRI reconstruction see \cite{Ribes_08}.\\
Standard regularization can improve the reconstructed image quality when the experimental conditions
are not too severe.
However, under a low magnetic field intensity with a high reduction factor, the problem becomes more difficult because of the low SNR during the acquisition process and poor estimation of coil sensitivity maps. In fact, in this paper, we are interested in a pMRI reconstruction problem under a low magnetic field with high reduction factor. We will show later that standard regularization techniques may not provide convincing solutions to this problem in such severe experimental conditions. Indeed, these methods assume Gaussianity of the solution in the spatial domain, which is not a quite good approximation since the empirical histogram of an MRI image voxels may be multimodal.
So, it is mandatory to design more sophisticated methods to reduce aliasing artifacts so as to ensure an acceptable reconstructed image quality. In our previous work \cite{chaari_08}, the appealing properties of wavelet transforms in generating sparse representations of regular images have been exploited to build a regularized reconstruction technique in the wavelet domain that reduces artifacts. This work is extended in this paper to further improve the reconstructed image quality, in particular by exploiting the artifact features. We will show that the choice of wavelet representations has also been motivated by the ability to employ tractable statistical models and the emergence of fast convex non differentiable optimization methods. 
The main contributions of this paper are the following:
\begin{itemize} 
\item the design of a novel wavelet-based inversion technique allowing us to incorporate useful constraints on the images to be reconstructed,
\item the proposition of an efficient convex optimization algorithm that minimizes the related regularized criterion.
\end{itemize}
This paper is organized as follows. In Section \ref{sec:background}, we give a brief overview of the SENSE method and its regularized version in the space domain. Section~\ref{sec:reg_WT}  describes the proposed wavelet-based regularization approach and the associated convex optimization algorithms. Results conducted on anatomical GRE and functional EPI data are provided and commented in Section~\ref{sec:simuls}. Finally, some conclusions and perspectives are drawn in Section~\ref{sec:concl}.
\section{Background} \label{sec:background}
Although pMRI is a 3D imaging technique, it proceeds with a slice by slice acquisition to get all the 3D volume, which simplifies the problem to a two-dimensional one. In what follows, we are only interested in reconstructing a given slice (2D image).
\subsection{Basic SENSE reconstruction }\hfill \\
An array of  $L$ coils  is employed to measure the spin density $\overline{\rho}$ into the object under investigation.\footnote{The overbar is used to distinguish the ``true'' data from a generic variable.} The signal $d_\ell$ received  by each coil $\ell$ ($1\leq \ell \leq L$) is the Fourier transform of the desired 2D field $\overline{\rho}$ weighted by the coil sensitivity profile $s_\ell$, evaluated at some locations $\vect{k}=(k_y,k_x)^{\trans}$ in
the $k$-space. The received signal $\widetilde{d}_\ell$ is therefore defined by the sampling scheme 
\begin{equation}
\widetilde{d}_\ell(\vect{k})=\int\overline{\rho}(\vect{r})s_\ell(\vect{r})e^{-\imath 2\pi
\vect{k}^{\trans}\vect{r}}\,d\vect{r} +\widetilde{n}(\vect{k}) \label{eq:signal}
\end{equation} 
where $\widetilde{n}(\vect{k})$ is a realization of an additive zero-mean Gaussian  noise and $\vect{r}=(y,x)^{\trans}$ is the spatial position in the image domain.
For the sake of simplicity, a Cartesian coordinate system is generally adopted in whole brain imaging.

In its simplest form, SENSE imaging amounts to solving a one-dimensional inversion problem due to the separability of the Fourier transform. Note however that it extends to a two-dimensional inversion problem for 3D imaging sequences like EVI (Echo Volume Imaging) \cite{Rabrait07} where undersampling occurs in two directions of the $k$-space.\\
Let $\Delta y=\frac{Y}{R}$ be the sampling period where $Y$ is the size of the FOV along the phase encoding direction, let $y$ be the position in the image domain along the same direction, $x$ the position in the image domain along the frequency encoding direction and $R \leq L$ the reduction factor. A 2D inverse Fourier transform allows us to recover the measured signal in the spatial domain. By accounting for the undersampling of the $k$-space by $R$, (\ref{eq:signal}) can be reexpressed in the following matrix form:
\begin{equation}
\vect{d}(\vect{r}) = \vect{S}(\vect{r}) \overline{\vect{\rho}}(\vect{r}) + \vect{n}(\vect{r})
\label{eq:matriciel}
\end{equation}
where 
\begin{equation}
\vect{S}(\vect{r})\,\eqdef\,\left(
\begin{array}{ccc}
s_1(y,x)&\ldots&s_1(y+(R-1)\Delta y,x)\\
\vdots&\vdots&\vdots\\
s_L(y,x)&\ldots&s_L(y+(R-1)\Delta y,x)\\
\end{array}
\right) \in \mathbb{C}^{L\times R}\nonumber,
\end{equation}

\begin{align}
\label{eq:defvrho}
& \overline{\vect{\rho}}(\vect{r})\,\eqdef\,
\left(
\begin{array}{c}
\overline{\rho}(y,x)\\
\overline{\rho}(y+\Delta y,x)\\
\vdots\\
\overline{\rho}(y+(R-1)\Delta y,x)\\
\end{array}
\right), 
 \quad 
\vect{d}(\vect{r})\,\eqdef\,
\left(
\begin{array}{c}
d_1(y,x)\\
d_2(y,x)\\
\vdots\\
d_L(y,x)\\
\end{array}
\right) \\ \nonumber
& \quad \mathrm{and}  \quad
\vect{n}(\vect{r})\,\eqdef\,\left(
\begin{array}{c}
n_1(y,x)\\
n_2(y,x)\\
\vdots\\
n_L(y,x)\\
\end{array}
\right).
\end{align}

In \eqref{eq:matriciel}, $(\vect{n}(\vect{r}))_{\vect{r}}$ is a sequence
of circular zero-mean Gaussian complex-valued random vectors. These noise vectors are i.i.d. and spatially independent with covariance matrix $\vect{\Psi}$ \cite{Pruessmann_K_99,Sijbers_07} of size $L \times L$. In practice, $\vect{\Psi}$ is estimated by acquiring data from all coils without radio frequency pulses, and its generic entry $\Psi(\ell_1,\ell_2)$ corresponding to the correlation between the two coils $\ell_1$ and $\ell_2$ is given by:
\begin{equation}
\Psi(\ell_1,\ell_2)=\frac{\sigma_n^2\sum_{(y,x)} s_{\ell_1}(y,x) s^{*}_{\ell_2}(y,x)}{\sqrt{(\sum_{(y,x)} |s_{\ell_1}(y,x)|^{2}) (\sum_{(y,x)} |s_{\ell_2}(y,x)|^{2})}}, \; \forall (\ell_1,\ell_2) \in \{1,\ldots,L\}^2
\end{equation} where $(\cdot)^{*}$ stands for the complex conjugate. 

Therefore, the reconstruction step consists in inverting (\ref{eq:matriciel}) and recovering $\overline{\vect{\rho}}(\vect{r})$ from $\vect{d}(\vect{r})$ at spatial positions $\vect{r}=(y,x)^{\trans}$.
Note that the data $(d_\ell)_{1\leq l \leq L}$ and the
unknown image $\overline{\rho}$ are complex-valued, although $|\overline{\rho}|$ is only considered for visualization purposes.

A simple reconstruction method also called the SENSE approach \cite{Pruessmann_K_99}, is based on the minimization of the Weighted Least Squares (WLS) criterion.
The objective is to find a vector $\widehat{\vect{\rho}}_{\rm WLS}(\vect{r})$ at each spatial location $\vect{r}$ such that:
\begin{align}
\widehat{\vect{\rho}}_{\rm WLS}(\vect{r})=&\arg\min_{\vect{\rho}(\vect{r})} \mathcal{J}_{\rm WLS}(\vect{\rho}(\vect{r})) \nonumber\\
=&\big(\vect{S}^{\hermit}(\vect{r})\vect{\Psi}^{-1}\vect{S}(\vect{r})\big)^{\sharp}\vect{S}^{\hermit}(\vect{r})\vect{\Psi}^{-1}\vect{d}(\vect{r})
    \label{eq:simple}
\end{align}
where $\mathcal{J}_{\rm WLS}(\vect{\rho}(\vect{r})) = \parallel\vect{d}(\vect{r})-\vect{S}(\vect{r})\vect{\rho}(\vect{r})\parallel^2_{\vect{\Psi}^{-1}}$, $(\cdot)^{\hermit}$ (resp. $(\cdot)^{\sharp}$) stands for the transposed complex conjugate (resp. pseudo-inverse)
and, $\|\cdot\|_{\vect{\Psi}^{-1}}= \sqrt{(\cdot)^{\hermit}\vect{\Psi}^{-1}(\cdot)}$ defines a norm on $\mathbb{C}^L$. \\
In practice, the performance of the SENSE method is limited because of the presence of i) distortions in the measurements $\vect{d}(\vect{r})$, ii) possible ill-conditioning of $\vect{S}(\vect{r})$ for some locations $\vect{r}$ and iii) the presence of errors in the estimation of $\vect{S}(\vect{r})$ mainly in the center and the brain air interface. These undesirable effects are illustrated in Fig.~\ref{fig:sense_tikh} (SENSE reconstruction) which shows aliasing artifacts in the reconstructed images for two values of the reduction factor $R=2$ and $R=4$. To increase the stability of the solution and reduce the artifacts, a Tikhonov regularization is usually applied.

\subsection{Tikhonov regularization}\hfill \\
As shown in \cite{Ying_L_04,Liang_02}, a Tikhonov regularization \cite{Tikhonov_63} can be applied to solve such an ill-posed inverse problem. 
The regularization process typically consists in computing $\widehat{\vect{\rho}}_{\rm PWLS}(\vect{r})$ as the minimizer of the following Penalized Weighted Least Squares (PWLS) criterion:
\begin{equation}
\mathcal{J}_{\rm PWLS}(\vect{\rho}(\vect{r}))=\mathcal{J}_{\rm WLS}(\vect{\rho}(\vect{r}))
+\kappa\parallel\!\vect{\rho}(\vect{r})-\vect{\rho}_\mathrm{r}(\vect{r})\!\parallel_{\vect{I}_R}^2.
\label{eq:Tk}
\end{equation}
where $\vect{I}_R$ is the $R$-dimensional identity matrix.
The  regularization parameter $\kappa >0$ ensures a balance between the closeness to the data and the penalty term, which
controls the deviation from a given reference vector $\vect{\rho}_\mathrm{r}(\vect{r})$. The solution $\widehat{\vect{\rho}}_{\rm PWLS}(\vect{r})$
admits the following closed-form expression:
\begin{align}
\begin{array}{ll}
\widehat{\vect{\rho}}_{\rm
PWLS}(\vect{r})=&\vect{\rho}_{\textrm{r}}(\vect{r})+ \\ 
&\big(\vect{S}^{\hermit}(\vect{r})\vect{\Psi}^{-1}\vect{S}(\vect{r})
+\kappa\vect{I}_R\big)^{-1}\vect{S}^{\hermit}(\vect{r})\vect{\Psi}^{-1}\big(\vect{d}(\vect{r})-\vect{S}(\vect{r})\vect{\rho}_{\mathrm{r}}(\vect{r})\big).
\end{array}
\label{eq:Tksol}
\end{align}
Note that the accuracy of the solution depends on the reference vector $\vect{\rho}_\mathrm{r}(\vect{r})$ and the choice of the regularization parameter $\kappa$. On the one hand, if $\kappa$ is close to zero, we tend to the WLS solution $\widehat{\vect{\rho}}_{\rm WLS}$. On the other hand, when $\kappa$ takes high values, the penalty term dominates and the solution tends to the reference $\vect{\rho}_{\mathrm{r}}$.
To select a relevant value for $\kappa$, some works have proposed an efficient choice of the regularization parameter \cite{Griesbaum_08,Engl_87,Rabrait07} by using for instance the discrepancy principle or the L-curve technique in order to enhance the reconstructed image quality. However, it is worth noticing in practice that any quadratic regularization procedure tends to introduce blurring effects, except for specific values of $\kappa$. An example of blurring effects can be seen in Fig~\ref{fig:sense_tikh} (Tikhonov reconstruction) which displays reconstructed full FOV images using Tikhonov regularization.
To overcome this limitation, one usually resorts to edge-preserving penalty terms \cite{wang_06,keeling_03,coulon_04}. Generally, these penalty terms are applied in the image domain and make the regularization more efficient by limiting blurring effects and preserving the image boundaries. However, the introduction of these terms may lead to a non-differentiable optimization problem which is not always easy to solve numerically.
Here, we propose to apply an edge-preserving regularization by proceeding in the Wavelet Transform (WT) domain \cite{Mallat_S_98}. We show that it is then possible to develop a fast and efficient algorithm to
solve the related convex non-differentiable optimization problem.

\section{Regularization in the WT domain}
\label{sec:reg_WT}
\subsection{Motivation}
When carefully analyzing SENSE-based reconstructed images (see Fig. \ref{fig:sense_tikh} for $R=4$), well spatially-localized artifacts appear as distorted curves with either very high or very low intensity. Consequently, we propose to look for an image representation where these localized transitions can be easily detected and hence attenuated. In this respect, the WT has been recognized as a powerful tool that enables a good space and frequency localization of useful information \cite{Mallat_S_98}. In the literature, many wavelet decompositions 
and extensions of these decompositions have been reported offering different features in order to provide sparse image representations. We can mention, for example, decompositions onto orthonormal dyadic wavelet bases \cite{Daubechies_88}
including the Haar transform \cite{Haar_10} as a special simple case, decompositions onto
biorthogonal dyadic wavelets \cite{cohen_92}, $M$-band wavelet representations \cite{steffen_93} and wavelet packet representations \cite{coifman_92}.

An appealing property of the resulting decomposition is that the statistical distribution  of the wavelet coefficients can be easily modelled in a realistic way. Hence, the Bayesian framework can be adopted to build an efficient reconstruction procedure.
Wavelet decompositions have been therefore investigated in a number of works in image processing like denoising \cite{Leporini_01, Muller_99, Heurta_05, daubechies_05} and deconvolution \cite{daubechies_05,Vonesch_08,chaux2_06} problems. In medical imaging, wavelet decompositions have also been investigated for instance in denoising \cite{Weaver_91,wang_06,Pizurica_06}, coil sensitivity map estimation \cite{lin_03} and encoding schemes \cite{gelman_96,Wendt_98} in MRI, activation detection in fMRI \cite{ruttimann_98,meyer_03,vandeville_04}, tissue characterization in ultrasound imaging \cite{mojsilovic_98} and tomographic reconstruction \cite{olson_92}. 
\subsection{Definitions and notations}
In what follows, $T$ stands for the WT operator. It corresponds to a discrete decomposition onto a separable 2D 
$M$-band wavelet basis performed over $j_{\mathrm{max}}$ resolution levels. 
The objective image $\overline{\rho}$ of size $Y \times X$ can be viewed as
an element of the Euclidean space $\mathbb{C}^K$ with $K = Y \times X$ endowed
with the standard inner product $\scal{\cdot}{\cdot}$ and norm $\|\cdot \|$. We recall that here we are interested in reconstructing one slice (2D image). For this reason, only 2D WT operators are investigated. In this context, the following notations are introduced.
\begin{definition} \hfill \\
Let $(e_k)_{1\le k \le K}$ be the considered discrete wavelet basis of the space $\mathbb{C}^{K}$. The wavelet decomposition operator $T$ is defined as the linear operator:
\begin{align}
T\colon\mathbb{C}^{K} &\to \mathbb{C}^{K}\\ \nonumber
\rho &\mapsto  (\scal{\rho}{e_k})_{1 \le k \le K}.
\end{align}
The adjoint operator $T^*$ serving for reconstruction purposes 
is then defined as the bijective linear operator:
\begin{align}
T^*\colon\mathbb{C}^{K}&\to \mathbb{C}^{K}\\ \nonumber
(\zeta_k)_{1 \le k\le K} &\mapsto \sum_{k=1}^K \zeta_k e_k.
\end{align}
\end{definition}
The resulting wavelet coefficient field of a target image function $\rho$ is defined by 
$\zeta =\big((\zeta_{a,k})_{1 \le k\le K_{j_\mathrm{max}}}, (\zeta_{o,j,k})_{1 \le j \le j_\mathrm{max}, 1\le k \le K_j}\big)$ 
where 
$K_{j}= KM^{-2j}$ is the number of wavelet coefficients 
in a given subband at resolution $j$ (by assuming that $Y$ and $X$ are multiple of $M^{j_{\mathrm{max}}}$) and the coefficients have been reindexed in such a way that
$\zeta_{a,k}$ denotes an approximation coefficient at resolution level $j_{\mathrm{max}}$ and $\zeta_{o,j,k}$ denotes a detail coefficient at resolution level $j$ and orientation $o \in \mathbb{O} =\{0,\ldots,M-1\}^2\setminus \{(0,0)\}$. In the dyadic case ($M=2$), there are three orientations corresponding to the
horizontal, vertical or diagonal directions.
Note that, when an orthonormal wavelet basis is considered, the adjoint operator
$T^*$ reduces to the inverse WT operator $T^{-1}$ and the operator norm $\|T\|$
of $T$ is equal to 1.

\subsection{Reconstruction procedure}

An estimate of the objective image $\overline{\rho}$ will be generated through
the reconstruction wavelet operator $T^*$. Let $\overline{\zeta}$ be the unknown
wavelet coefficients such that $\overline{\rho} = T^* \overline{\zeta}$.
We aim at building  an estimate $\widehat{\zeta}$ of the vector of coefficients
$\overline{\zeta}$ from the observations $\vect{d}$. To this end, we derive a Bayesian approach relying on suitable priors on the wavelet coefficients. 

\subsubsection{Likelihood}\label{sec:likelihood}\hfill \\
Since the noise has been assumed i.i.d. circular Gaussian with zero-mean and between-coil correlation matrix $\vect{\Psi}$, the data likelihood can be expressed as:
\begin{equation}
g(\vect{d} \mid T^*\zeta)=\prod_{\mathbf{r}\in \{1,\ldots,Y\}\times \{1,\ldots,X\}}g(\vect{d}(\mathbf{r}) \mid \vect{\rho}(\mathbf{r})) \propto
\exp\left(-\mathcal{J}_{\rm L}(\rho)\right)\label{eq:likelihood}
\end{equation}
where $\rho = T^*\zeta$ and
\begin{equation}
\mathcal{J}_{\rm L}(\rho) =\sum_{\mathbf{r}\in \{1,\ldots,Y\}\times \{1,\ldots,X\}} \mathcal{J}_{\rm WLS}(\vect{\rho}(\mathbf{r})) \nonumber
\label{eq:defG}
\end{equation}

\subsubsection{Prior}\label{sec:prior}\hfill \\
Let us denote by $f$ the prior probability density function (PDF) of the image in the wavelet domain. By analyzing the correlation between the real and imaginary parts of the wavelet coefficients, very low values of the correlation factors have been found for medical pMRI images. Table \ref{tab:corr} illustrates the values of the correlation coefficient.
\begin{table}[!ht]
\centering
\caption{Correlation between real and imaginary parts of wavelet coefficients over 3 resolution levels.}
\vspace*{0.3cm}
\begin{tabular}{|p{1.5 cm}|p{1.5 cm}|p{3 cm}|p{2 cm}|}
  \hline
\multicolumn{2}{|c|}{}& Approximation & Detail \\
  \hline 
&$j=1$ & -0.004  & -0.139 \\
  \cline{2-4}
  Slice 5 &$j=2$ & -0.031  & -0.140 \\
  \cline{2-4}
  &$j=3$ & -0.026 & -0.153 \\
  \hline
&$j=1$ & -0.111  & -0.077 \\
  \cline{2-4}
  Slice 9 &$j=2$ & -0.117  & -0.032 \\
  \cline{2-4}
  &$j=3$ & -0.122 & -0.031 \\
  \hline
\end{tabular}
\label{tab:corr}
\end{table}

\begin{center}
\begin{figure}[!ht]
\centering
\includegraphics[width=10cm, height=4cm]{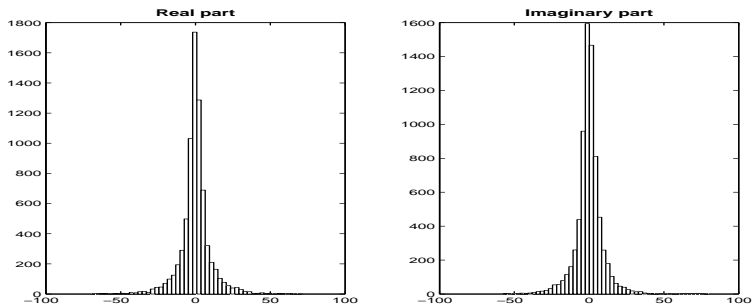}
\caption{Example of empirical histogram of wavelet coefficients.}
\label{fig:hist}
\end{figure}
\end{center}

These observations have motivated the choice of independent priors for the real and imaginary parts of the wavelet coefficients, the marginal distributions of which can be separately studied. For the sake of simplicity, we also assume that the real (resp. imaginary) parts of the wavelet coefficients are i.i.d. in each subband. Their statistical characteristics may however vary between two distinct subbands.
Furthermore, by looking at the empirical distributions of the real and imaginary parts of the considered wavelet coefficients, we have noticed that their empirical histograms are well-fitted by a Generalized Gauss-Laplace distribution. The histograms present a single mode and their shape vary between the Gaussian and Laplacian ones (see Fig.~\ref{fig:hist}). The corresponding PDF is:
\begin{equation}
\label{eq:loi_x2}
\forall \xi \in \mathbb{R}, \quad f(\xi;\alpha,\beta)= \sqrt{\frac{\beta}{2\pi}}
\dfrac{e^{-(\alpha|\xi|+\frac{\beta}{2} \xi^2+\frac{\alpha^2}{2\beta})}}{\mathrm{erfc}(\frac{\alpha}{\sqrt{2\beta}})},
\end{equation} 
where $\alpha\in \RR_+$ and $\beta$ in $\RR_+^*$ are hyper-parameters to be estimated.

At the coarsest resolution level $j_{\mathrm{max}}$, the distributions of both the real and imaginary parts of the approximation coefficients can be modelled by a Gaussian distribution since they belong to a low frequency subband.

\subsubsection{Bayesian inference}\label{sec:bayesian}\hfill \\
Based on the prior and the likelihood given hereabove, the MAP estimator is  computed by minimizing the following criterion:
\begin{align}
\widehat{\zeta} &= \arg\max_{\zeta\in\mathbb{C}^K} \big(\ln f(\zeta)+\ln g(\vect{d} \mid T^{*} \zeta)\big) \nonumber  \nonumber \\ 
&= \arg\min_{\zeta\in\mathbb{C}^K}\mathcal{J}_{\rm WT}(\zeta) \label{eq:argmax}
\end{align}
where
$\mathcal{J}_{\rm WT}(\zeta)=\mathcal{J}_{L}(T^*\zeta) + \mathcal{J}_{P}(\zeta)$ and
\begin{equation}
\mathcal{J}_P(\zeta) = \sum_{k=1}^{K_{j_\mathrm{max}}} \Phi_{a}(\zeta_{a,k}) +\sum_{o\in \mathbb{O}}\sum_{j=1}^{j_{\mathrm{max}}} \sum_{k=1}^{K_j} 
\Phi_{o,j}(\zeta_{o,j,k}) \nonumber
\label{eq:J20}
\end{equation}
with
\begin{align}
\Phi_{a}(\zeta_{a,k}) &=
\frac{(\mathrm{Re}(\zeta_{a,k})-\mu^{\mathrm{Re}})^{2}}{2\sigma_{\mathrm{Re}}^{2}}+
\frac{(\mathrm{Im}(\zeta_{a,k})-\mu^{\mathrm{Im}})^{2}}{2\sigma _{\mathrm{Im}}^{2}},\label{eq:J21}\\ \nonumber
\Phi_{o,j}(\zeta_{o,j,k}) &=
\alpha_{j,o}^{\mathrm{Re}}|\mathrm{Re}(\zeta_{o,j,k})| +
\frac{\beta_{j,o}^{\mathrm{Re}}}{2}|\mathrm{Re}(\zeta_{o,j,k})|^2
+\alpha_{j,o}^{\mathrm{Im}}|\mathrm{Im}(\zeta_{o,j,k})| +
\frac{\beta_{j,o}^{\mathrm{Im}}}{2}|\mathrm{Im}(\zeta_{o,j,k})|^2,\label{eq:J22}
\end{align}
$(\alpha_{j,o}^{\mathrm{Re}},\alpha_{j,o}^{\mathrm{Im}})\in(\mathbb{R}_+^*)^2$, $(\beta_{j,o}^{\mathrm{Re}},\beta_{j,o}^{\mathrm{Im}})\in(\mathbb{R}_+^*)^2$, $(\mu^{\mathrm{Re}},\mu^{\mathrm{Im}})\in \mathbb{R}^2$ and $(\sigma_{\mathrm{Re}},\sigma_{\mathrm{Im}}) \in \mathbb{R}_+^2$.

Hereabove, $\mathrm{Re}(\cdot)$ and $\mathrm{Im}(\cdot)$ (or $\cdot^{\mathrm{Re}}$ and
$\cdot^{\mathrm{Im}}$) stand for the real and imaginary parts, respectively. 
It is clear that $\mathcal{J}_{\rm WT}$ is convex.
However, the optimization procedure cannot rely on conventional convex optimization techniques like the pseudo-conjugate gradient: although $\mathcal{J}_L$ is differentiable with a Lipschitz-continuous gradient, $\mathcal{J}_P$ is not differentiable. This is a main difficulty which is frequently encountered in inverse problems involving sparsity promoting priors \cite{harikumar_96,tropp_06}.  Therefore, we propose to apply a generalized form of the iterative optimization procedure developed in \cite{Chaux_C_07} which extends the one in \cite{Daub_04} and is based on the forward-backward algorithm.

\subsection{Forward-Backward optimization algorithm}
\label{sec:algo} 
The minimization of $\mathcal{J}_{\rm WT}$ is performed by resorting to the concept of proximity operators \cite{Moreau_65}, which was found to be fruitful
in a number of recent works in convex optimization \cite{Chaux_C_07,Combettes_05}.
In what follows, we give a brief overview of this key tool 
for solving our optimization problem. 
\begin{definition} {\rm \cite{Moreau_65}} \label{def:prox}

Let $\Gamma_0(\chi)$ be the class of lower semicontinuous convex functions from a separable real Hilbert space $\chi$ to $]-\infty,+\infty]$ and let $\varphi \in \Gamma_0(\chi)$. For every $\mathsf{x} \in \chi$, the function $\varphi(\cdot)+\Vert \cdot-\mathsf{x} \Vert^2/2$ achieves its infimum at a unique point denoted by $\mathrm{prox}_{\varphi}\mathsf{x}$. The operator $\mathrm{prox}_{\varphi}\; : \; \chi \rightarrow \chi$ is the proximity operator of $\varphi$.
\end{definition}
Here is an example of proximity operators which will be used in our approach.
\begin{example}\hfill \\
Consider the following function:
\begin{align}
\varphi \colon \mathbb{R}& \to \mathbb{R}\\ \nonumber
\xi &\mapsto \alpha|\xi-\mu|+\frac{\beta}{2} (\xi-\mu)^2
\end{align}
with $\alpha \in \RR_+$, $\beta \in \RR_+^*$ and $\mu \in \mathbb{R}$. The associated proximity operator is given by:
\begin{equation}
\mathrm{prox}_{\varphi} \xi=  \dfrac{\mathrm{sign}(\xi-\mu)}{\beta+1}\max\{|\xi-\mu|- \alpha,0\}+\mu
\end{equation}
where the $\mathrm{sign}$ function is defined as follows:
\begin{equation}
\forall \xi\in\mathbb{R},\qquad  \mathrm{sign}(\xi)= \begin{cases} +1 & \text{if} \; \xi \geq 0\\
-1 & \text{otherwise.} \end{cases} \nonumber
\end{equation}
\end{example}

In this work, as the observed data are complex-valued, we generalize the definition of proximity operators to a class of convex functions defined for complex-valued variables.\\ For the function
\begin{align}
\Phi \colon \mathbb{C}^K &\to ]-\infty,+\infty]\\ \nonumber
x &\mapsto \phi^{\mathrm{Re}}(\mathrm{Re}(x))+ \phi^{\mathrm{Im}}(\mathrm{Im}(x)),
\end{align} 
where $\phi^{\mathrm{Re}}$ and $\phi^{\mathrm{Im}}$ are functions in $\Gamma_0(\RR^K)$ and $\mathrm{Re}(x)$ (resp. $\mathrm{Im}(x)$) is the vector
of the real parts (resp. imaginary parts) of the components of $x\in \mathbb{C}^K$, the proximity operator is defined as:
\begin{align}
\mathrm{prox}_{\Phi} \colon \mathbb{C}^K & \to \mathbb{C}^K \\ \nonumber
x &\mapsto \mathrm{prox}_{\phi^{\mathrm{Re}}}(\mathrm{Re}(x))+\imath \mathrm{prox}_{\phi^{\mathrm{Im}}}(\mathrm{Im}(x)).
\label{eq:defproxc}
\end{align}
Here is an example of proximity operator for a function of a complex-valued variable.
\begin{example}\label{ex:pcomp}\hfill \\
Consider the following function:
\begin{align}
\Phi \;:\;\mathbb{C} &\to \mathbb{R}\\ \nonumber
\xi & \mapsto \alpha^{\mathrm{Re}}|\mathrm{Re}(\xi-\mu)|+\frac{\beta^{\mathrm{Re}}}{2}\big(\mathrm{Re}(\xi-\mu)\big)^2 \\ \nonumber
&+ \alpha^{\mathrm{Im}}|\mathrm{Im}(\xi-\mu)|+
\frac{\beta^{\mathrm{Im}}}{2}|\mathrm{Im}(\xi-\mu)|^2
\end{align}
with $(\alpha^{\mathrm{Re}},\alpha^{\mathrm{Im}})\in (\RR_+)^2$, $(\beta^{\mathrm{Re}},\beta^{\mathrm{Im}}) \in (\mathbb{R}_+^*)^2$ and $\mu \in \mathbb{C}$. The associated proximity operator is:
\begin{multline}
\mathrm{prox}_{\Phi} \xi=  \dfrac{\mathrm{sign}(\mathrm{Re}(\xi-\mu))}{\beta^{\mathrm{Re}}+1}\max\{|\mathrm{Re}(\xi-\mu)|- \alpha^{\mathrm{Re}},0\}\\ + \imath \dfrac{\mathrm{sign}(\mathrm{Im}(\xi-\mu))}{\beta^{\mathrm{Im}}+1}\max\{|\mathrm{Im}(\xi-\mu)|- \alpha^{\mathrm{Im}},0\}+\mu.
\end{multline}
\end{example}

Based on these definitions, and by extending the algorithm in \cite{Chaux_C_07} to the complex case, a minimizer of $\mathcal{J}_{\rm WT}$ can then be iteratively computed according to Algorithm \ref{algo:wt}.
Note that in this algorithm, the expressions of $\mathrm{prox}_{\gamma _{n} \Phi_{a}}$ and $\mathrm{prox}_{\gamma _{n} \Phi_{o,j}}$ at each iteration $n$ are provided by Example~\ref{ex:pcomp}.
It can also be noticed that $\lambda_n$ and $\gamma_n$ respectively correspond to relaxation and step-size parameters. 

\begin{algorithm}
Let $(\gamma_n)_{n>0}$ and $(\lambda_n)_{n>0}$ be sequences of positive reals.
\caption{2D-slice wavelet-based regularized reconstruction}
\label{algo:wt}
\begin{algorithmic}[1]
\STATE Initialize $\zeta^{(1)}$. Set $n=1$, $\epsilon \in (0,1)$
and $\mathcal{J}^{(0)}=0$.
\REPEAT
	\STATE\label{s:a1-1} Reconstruct the image by setting $\rho^{(n)} = T^* \zeta^{(n)}$.
	\STATE Calculate the image $u^{(n)}$ such that:\\
	$\forall \mathbf{r}\in \{1,\ldots,Y\}\times \{1,\ldots,X\}$,\\ $\vect{u}^{(n)}(\mathbf{r})=2\vect{S}^{\hermit}(\mathbf{r})\vect{\Psi}^{-1}\left(\vect{S}(\mathbf{r})
\boldsymbol{\rho}^{(n)}(\mathbf{r})-\vect{d}(\mathbf{r})\right)$,\\
where the vector $\vect{u}^{(n)}(\mathbf{r})$ is defined from $u^{(n)}$ in the same way
as $\overline{\boldsymbol{\rho}}(\mathbf{r})$ is defined from $\overline{\rho}$ (see \eqref{eq:defvrho}).
	\STATE Determine the wavelet coefficients $\upsilon^{(n)} = T u^{(n)}$ of $u^{(n)}$.
	\STATE Update the approximation coefficients of the reconstructed image:\\
$\forall k \in \KK_{j_\mathrm{max}}$,
$\zeta_{a,k}^{(n+1)} = \zeta_{a,k}^{(n)}+\lambda_{n}\Big(
\mathrm{prox}_{\gamma _{n} \Phi_{a}}
(\zeta_{a,k}^{(n)}-\gamma_{n} \upsilon^{(n)}_{a,k})-\zeta_{a,k}^{(n)}\Big)$.
	\STATE Update the detail coefficients of the reconstructed image:
\begin{multline*}
\forall o \in \mathbb{O},
\forall j \in \{1,\ldots,j_{\mathrm{max}}\}, \forall k \in \KK_j,\\
\zeta_{o,j,k}^{(n+1)} = \zeta_{o,j,k}^{(n)}+\lambda_{n}\Big(
\mathrm{prox}_{\gamma _{n} \Phi_{o,j}}
(\zeta_{o,j,k}^{(n)}-\gamma_{n}\upsilon^{(n)}_{o,j,k})-\zeta_{o,j,k}^{(n)}\Big) .
\end{multline*}
\STATE  Compute  $\mathcal{J}^{(n)}  = \mathcal{J}_{\rm WT}(\zeta^{(n)})$.
\STATE \label{s:a1-f} $n \leftarrow n+1$
\UNTIL {$| \mathcal{J}^{(n-1)}-\mathcal{J}^{(n-2)}| \le \varepsilon \mathcal{J}^{(n-1)}$}
\RETURN $\rho^{(n)}=T^*\zeta^{(n)}$
\end{algorithmic}
\end{algorithm}

For 3D volume reconstruction we iterate over slices, and for functional data, we then iterate over volumes independently (see discussion in Section~\ref{sec:concl}).  
\subsection{Convergence of Algorithm~\ref{algo:wt}}
For every $\mathbf{r}\in \{1,\ldots,Y\}\times 
\{1,\ldots,X\}$, let $\theta_{\mathbf{r}}\ge 0$ be the maximum eigenvalue of the Hermitian positive semi-definite matrix $\vect{S}^{\hermit}(\mathbf{r})\vect{\Psi}^{-1}\vect{S}(\mathbf{r})$ and let $\theta = \max_{\mathbf{r}\in \{1,\ldots,Y\}\times 
\{1,\ldots,X\}}\theta_{\mathbf{r}} > 0$. To guarantee the convergence of Algorithm \ref{algo:wt}, the step-size and relaxation parameters are subject to the following conditions:

\begin{assumption}\  \label{as:step}
\begin{enumerate}
\item $\inf_{n>0}\gamma_n > 0$ and  $\sup_{n>0}\gamma_n< \frac{1}{\theta\|T\|^2}$,
\item $\inf_{n>0} \lambda_n > 0$ and $\sup_{n>0} \lambda_n \le 1$.
\label{ass:cv}
\end{enumerate}
\end{assumption}
More precisely, the following result can be shown:

\begin{proposition}\label{prop:cvl}
Under Assumption \ref{as:step}, the sequence $(\zeta^{(n)})_{n>0}$
built when iterating Steps \ref{s:a1-1} to \ref{s:a1-f} of Algorithm \ref{algo:wt} converges linearly to the unique
solution $\widehat{\zeta}$  to Problem \eqref{eq:argmax}.
\end{proposition}
\textit{Proof:} See Appendix A.\\

Results obtained with the proposed iterative method are provided and discussed in Section \ref{sec:simuls}. In these results, it can be noticed that some artifacts remain in the reconstructed full FOV images (see the image obtained with Algorithm 1 in Fig.~\ref{fig:algos12}) because they present extremely large volumes. In the next section, we therefore present an extension based on additional constraints that enable a very accurate reconstruction and cancellation of such artifacts.

\subsection{Constrained wavelet-based regularization}
\label{sec:constrained} 

We propose to extend our approach
by incorporating an additional constraint in the method described hereabove in order to better regularize artifact regions. 

\subsubsection{New optimality criterion}
As artifacts appear as curves with very high or very low intensity, we propose to set local lower and upper bounds on the image intensity values in artifact areas. 
These bounds define the nonempty closed convex set:
\begin{equation}
C = \{\rho \in \mathbb{C}^K \mid \forall \vect{r} \in \{1,\ldots,Y\}\times
\{1,\ldots,X\}, \rho(\vect{r}) \in C_\vect{r} \}
\end{equation}
where the constraint introduced on the range values at position $\vect{r}\in \{1,\ldots,Y\}\times \{1,\ldots,X\}$ is modelled by
\begin{equation}\label{eq:Cr}
C_\vect{r}=\{\xi\in\mathbb{C}\mid \mathrm{Re}(\xi)\in\mathbb{I}^{\mathrm{Re}}_\vect{r}, \mathrm{Im}(\xi)\in\mathbb{I}^{\mathrm{Im}}_\vect{r}\}, 
\end{equation}
with $\mathbb{I}^{\mathrm{Re}}_\vect{r}=[\mathrm{I}^{\mathrm{Re}}_{\mathrm{min},\vect{r}},\mathrm{I}^{\mathrm{Re}}_{\mathrm{max},\vect{r}}]$
and $\mathbb{I}^{\mathrm{Im}}_\vect{r}=[\mathrm{I}^{\mathrm{Im}}_{\mathrm{min},\vect{r}},\mathrm{I}^{\mathrm{Im}}_{\mathrm{max},\vect{r}}]$.

When taking into account the additional constraint, the constrained criterion in \eqref{eq:argmax} becomes:
\begin{equation}
 \mathcal{J}_{\rm CWT}(\zeta)= \mathcal{J}_L(\zeta)+\mathcal{J}_P(\zeta) + i_{C^*}(\zeta),
\end{equation}
where 
\begin{equation}
C^* = \{\zeta \in \mathbb{C}^K \mid T^*\zeta \in C\} \nonumber
\end{equation}
and $i_{C^*}$ is the indicator function of the closed convex set $C^*$
defined by
\begin{equation}
\forall \zeta\in\mathbb{C}^K,\qquad i_{C^*}(\zeta)= \begin{cases} 0 & 
\mbox{if $\zeta \in C^*$}\nonumber \\
+ \infty & \text{otherwise.} \end{cases}
\end{equation}
The minimization problem can then be rewritten as:
\begin{equation}
\label{eq:constpb}
\mathrm{Find} \; \min_{\zeta\in\mathbb{C}^K } \mathcal{J}_L(\zeta)+\mathcal{J}_P(\zeta) + i_{C^*}(\zeta).
\end{equation}

\subsubsection{Computation of the new estimator}

Conceptually, in order to solve the minimization problem in \eqref{eq:constpb}, the forward-backward iteration in \eqref{eq:FW} (in Appendix~A) has to be updated according to:
\begin{equation}
\zeta^{(n+1)} = \zeta^{(n)}+\lambda_{n}\Big(
\mathrm{prox}_{\gamma _{n} \mathcal{J}_P+i_{C^*}}
\big(\zeta^{(n)}-\gamma_{n}\nabla \mathcal{J}_L(\zeta^{(n)})\big)-\zeta^{(n)}\Big).
\end{equation}
The main difficulty here is that the proximity operator of
$\gamma _{n} \mathcal{J}_P+i_{C^*}$ does not have a closed form.
However, from the definition of the proximity operator, we get
\begin{equation}
\forall \zeta \in \mathbb{C}^K,\qquad
\mathrm{prox}_{\gamma _{n} \mathcal{J}_P+i_{C^*}}(\zeta)
= \arg\min_{\zeta'\in\mathbb{C}^K} \gamma _{n} \mathcal{J}_P(\zeta')+\mathcal{J'}_{\zeta}(\zeta')
\label{eq:minprox}
\end{equation}
where
\begin{equation}
\mathcal{J'}_{\zeta} (\cdot)= \frac{1}{2} \|\cdot-\zeta\|^2+i_{C^*}(\cdot).
\end{equation}
Although $\mathrm{prox}_{\gamma _{n} \mathcal{J}_P+i_{C^*}}$ does not take a simple expression, the proximity operator of $\gamma _{n} \mathcal{J}_P$ is expressed by \eqref{eq:proxJ2} (in Appendix~A)
and the proximity operator of $\mathcal{J'}_{\zeta}$ is easily determined.
Indeed, from simple calculations, we have
\begin{equation}
\forall \zeta' \in \mathbb{C}^K,\qquad
\prox_{\mathcal{J'}_{\zeta}} (\zeta') = P_{C^*}\left(\frac{\zeta'+\zeta}{2}\right)
\end{equation}
where $P_{C^*}$ is the projection onto the convex set $C^*$.
In turn, provided that the considered wavelet basis is orthonormal, the projection onto $C^*$ of $\zeta'\in \mathbb{C}^K$ is obtained by performing the wavelet decomposition of the projection of 
$\rho' = T^*\zeta'$ onto $C$. The latter projection is 
\begin{equation}
P_C(\rho') = \big(P_{C_{\vect{r}}}(\rho'(\vect{r})\big)_{\vect{r} \in \{1,\ldots,Y\}\times \{1,\ldots,X\}}
\end{equation}
where, for every $\vect{r} \in \{1,\ldots,Y\}\times \{1,\ldots,X\}$,
\begin{equation}
\forall \xi \in \mathbb{C},\qquad
\mathrm{Re}\big(P_{C_{\vect{r}}}(\xi)\big) =\begin{cases}
\mathrm{I}^{\mathrm{Re}}_{\mathrm{min},\vect{r}} & \mbox{if $\mathrm{Re}(\xi) < \mathrm{I}^{\mathrm{Re}}_{\mathrm{min},\vect{r}}$}\\
\mathrm{I}^{\mathrm{Re}}_{\mathrm{max},\vect{r}} & \mbox{if $\mathrm{Re}(\xi) > \mathrm{I}^{\mathrm{Re}}_{\mathrm{max},\vect{r}}$}\\
\xi & \mbox{otherwise,}\\
\end{cases}
\label{eq:ressep2}
\end{equation}
a similar expression being used to calculate $\mathrm{Im}\big(P_{C_{\vect{r}}}(\xi)\big)$.

Knowing $\prox_{\gamma_n \mathcal{J}_P}$ and $\prox_{\mathcal{J'}_{\zeta}}$,
$\mathrm{prox}_{\gamma _{n} \mathcal{J}_P+i_{C^*}}\zeta$ can be computed in an iterative manner by solving the optimization problem in \eqref{eq:minprox} with the
Douglas-Rachford algorithm  \cite{pustelnik_08,Combettes_07}.
More precisely, we have:
\begin{proposition}\hfill \\
\label{prop:cv}
Set $\eta^{(0)} \in \mathbb{C}^K$ and construct for all $m \in \NN$:
\begin{equation}
\begin{cases}	
\eta^{(m+\frac{1}{2})} = \prox_{\mathcal{J'}_{\zeta}} \eta^{(m)}  \\
\eta^{(m+1)} = \eta^{(m)} + \tau\big(\prox_{\gamma_n \mathcal{J}_P}(2 \eta^{(m+\frac{1}{2})} - \eta^{(m)}) -\eta^{(m+\frac{1}{2})}\big),
\end{cases}
\label{eq:DR}
\end{equation}
where $\tau \in ]0,2[$.
Then, $(\eta^{(m+\frac{1}{2})})_{m\in \NN}$ converges to $\prox_{\gamma_n\mathcal{J}_{P}+i_{C^*}} \zeta$.
\end{proposition}
Inserting this extra iterative step in the forward-backward algorithm and using the expressions of $\prox_{\gamma_n \mathcal{J}_P}$ and $\prox_{\mathcal{J'}_{\zeta}}$ lead to Algorithm~\ref{algo:main}.\\
At iteration $n$, $M_n$ is the number of times the Douglas-Rachford step is run. According to the convergence analysis conducted in \cite[Prop. 4.2]{pustelnik_08}, if $M_n$ is chosen large enough and Assumption \ref{as:step} holds, 
iterating Steps \ref{s:a2-1} to \ref{s:a2-f} of the algorithm guarantees the  convergence to the unique solution to Problem \eqref{eq:constpb}. Note however that 
\cite[Prop. 4.2]{pustelnik_08} does not provide a practical guideline for the choice of $M_n$. The practical rule we chose is explained in Section \ref{sec:anat_data}.\\
The improvements resulting from the use of this constrained approach are illustrated in the next section.
\newpage
\begin{algorithm}[!ht]
Let $(\gamma_n)_{n>0}$ and $(\lambda_n)_{n>0}$ be sequences of positive reals, let $(M_n)_{n>0}$ be a sequence of positive integers and set $\tau\in ]0,2]$.
\caption{Constrained 2D-slice wavelet-based regularized reconstruction}
\label{algo:main}
\begin{algorithmic}[1]
\STATE Initialize $\zeta^{(1)}$. Set $n=1$, $\varepsilon \in (0,1)$
and $\mathcal{J}^{(0)} = 0$.
\REPEAT
	\STATE Same as for Algorithm 1 \label{s:a2-1} 
	\STATE Same as for Algorithm 1
	\STATE Same as for Algorithm 1
	\STATE Initialize the Douglas-Rachford algorithm by setting\\ $\eta^{(n,0)}=\zeta^{(n)}-\gamma_n \upsilon^{(n)}$.\\
	\STATE Douglas-Rachford iterations:\\
	\begin{itemize} 
		\FOR {$m=0$ to $M_n-1$ }
		\STATE Compute $\displaystyle \eta^{(n,m+\frac{1}{2})} = P_{C^*}\Big(\frac{\eta^{(n,m)} +\zeta^{(n)}}{2}\Big)$; \\
		\STATE Update the approximation components of $\eta^{(n,m)}$:\\
		$\forall k \in \KK_{j_\mathrm{max}}$, $\eta^{(n,m+1)}_{a,k} = \eta^{(n,m)}_{a,k} + \tau\big(\prox_{\gamma_n \Phi_a}(2 \eta^{(n,m+\frac{1}{2})}_{a,k} - \eta^{(n,m)}_{a,k}) -\eta^{(n,m+\frac{1}{2})}_{a,k}\big)$,
		\STATE Update the detail components of $\eta^{(n,m)}$:
		\vspace*{-0.5cm}
		\begin{multline*}
		\forall o \in \mathbb{O},
		\forall j \in \{1,\ldots,j_{\mathrm{max}}\}, \forall k \in \KK_j,\\
		\eta^{(n,m+1)}_{o,j,k} = \eta^{(n,m)}_{o,j,k} + \tau\big(\prox_{\gamma_n \Phi_{o,j}}(2 \eta^{(n,m+\frac{1}{2})}_{o,j,k} - \eta^{(n,m)}_{o,j,k}) -\eta^{(n,m+\frac{1}{2})}_{o,j,k}\big);
		\end{multline*}			
		\STATE If $\eta^{(n,m+1)} = \eta^{(n,m)}$, goto \ref{s:suiv}.
		\ENDFOR
	\end{itemize}
	\STATE  \label{s:suiv}Update the wavelet coefficients of the reconstructed image:\\
$\zeta^{(n+1)} = \zeta^{(n)} + \lambda_n (\eta^{(n,m+\frac{1}{2})}-\zeta^{(n)})$.
\STATE Compute $\mathcal{J}^{(n)}=\mathcal{J}_{\rm CWT}(\zeta^{(n)})$.
	\STATE \label{s:a2-f} $n \leftarrow n+1$.
\UNTIL {$| \mathcal{J}^{(n-1)}-\mathcal{J}^{(n-2)}| \le \varepsilon\mathcal{J}^{(n-1)}$}
\RETURN $\rho^{(n)}=T^*\zeta^{(n)}$
\end{algorithmic}
\end{algorithm}

\newpage

\section{Experimental results}
\label{sec:simuls}
Experiments have been conducted on real data sets containing $256\times256\times14$ Gradient-Echo anatomical and $64\times64\times30$ EPI functional images with respectively $0.93\times0.93\times8$~(mm$^3$) and $3.75\times3.75\times3$~(mm$^3$) spatial resolution. 
For Gradient-Echo data, we have TE/TR = $10/500$~(ms) with BW=$31.25$~(kHz), while for EPI data we have TE/TR = $30/2400$~(ms) and BW=$125$~(kHz).
These images have been acquired using reduction factors $R=2$ and $R=4$ to show the effectiveness of our approach on image reconstruction after a significant reduction of the global imaging time. The employed scanner was Signa 1.5 Tesla GE Healthcare with an eight-channel head coil.
The standard deviation $\sigma_n$ of the acquisition noise was estimated based on the background region of the measured images where there is no signal corresponding to the region of interest (in our case the head). In order to proceed to an objective assessment of the reconstruction accuracy, a full FOV image $\rho_{\rm ref}$ is available through a non accelerated acquisition for $R=1$.
\subsection{Anatomical data}\label{sec:anat_data}
\begin{figure}[!ht]
\centering
\includegraphics[width=8cm, height=4.4cm]{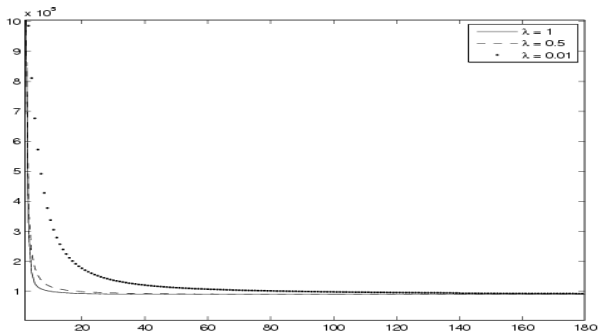}
\caption{Convergence speed of the optimization algorithm w.r.t. the choice of the relaxation parameter $\lambda$ for $j_{\mathrm{max}}=3$.}
\label{fig:critere1}
\end{figure}

In this experiment, dyadic ($M = 2$) Symmlet orthonormal wavelet bases \cite{daubechies_92} associated with filters of length 8 have been used over $j_{\mathrm{max}}=3$ resolution levels. For the wavelet coefficients, the priors described in Section \ref{sec:bayesian} have been employed. The related hyper-parameters have been estimated based on a reference image $\rho_{\rm ref}$ using the maximum likelihood estimator. A couple of hyper-parameters $(\alpha,\beta)$ is estimated for real/imaginary parts of each subband.
Algorithm~\ref{algo:wt} was then used to reconstruct full FOV images. For the sake of simplicity, constant values  along the algorithm iterations have been adopted for the relaxation parameter $\lambda_n$ and the step-size $\gamma_n$. It has been noticed experimentally that $\lambda_n \equiv 1$ is the optimal value of the relaxation parameter in terms of convergence rate (see Fig.~\ref{fig:critere1}). 

A value of the step-size closed to the allowed maximum value in  Assumption \ref{as:step} has also been observed to be beneficial to the convergence rate. After computing the 
constant $\theta$ related to the considered sensivity map, $\gamma_n$ was 
thus chosen equal to $1.99/\theta = 12.83$. 

\begin{figure}
\centering
\vspace*{0.1cm}
\begin{tabular}{c c c}
& $R=2$&$R=4$\\
 & \includegraphics[width=4cm, height=4cm]{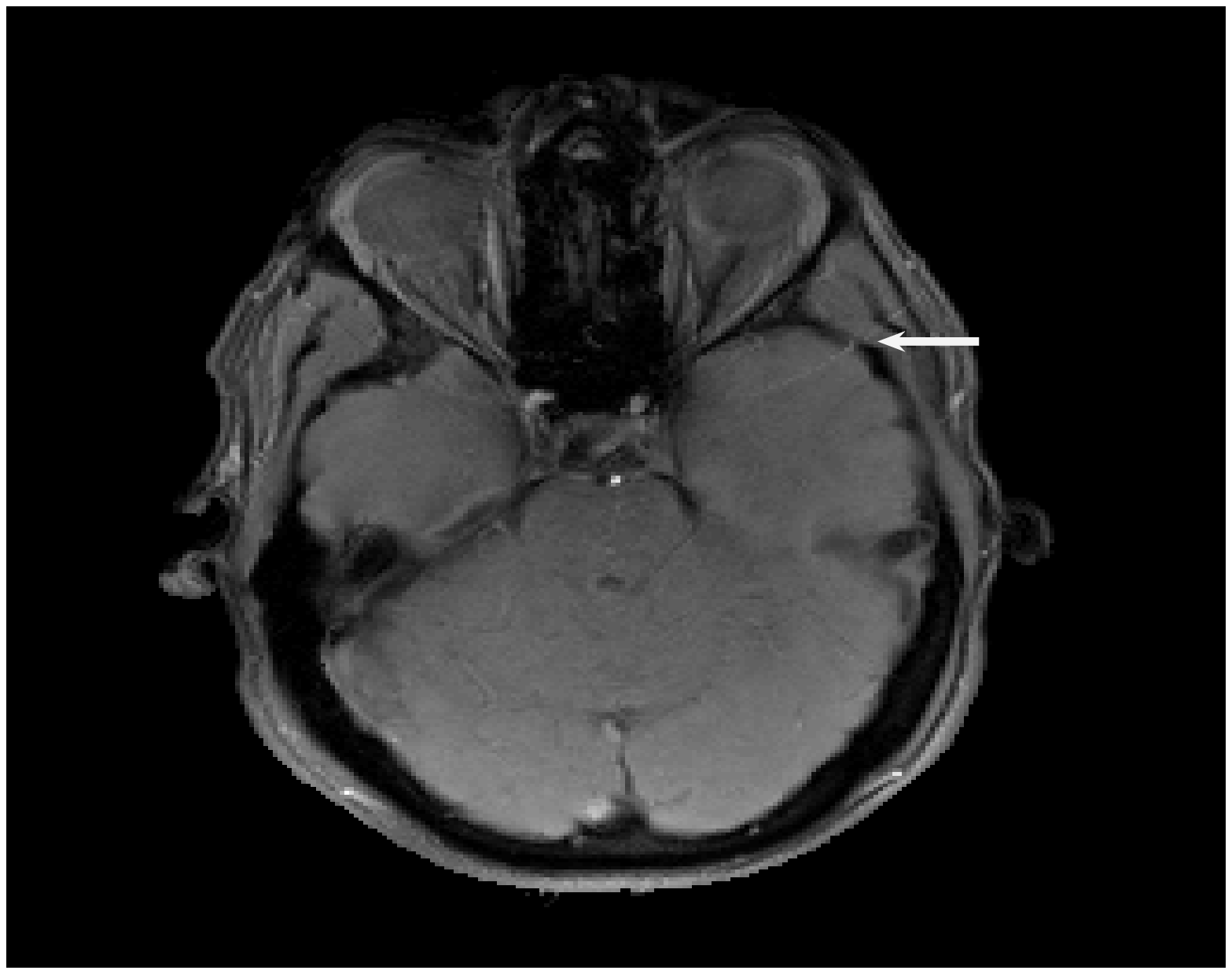} & \includegraphics[width=4cm, height=4cm]{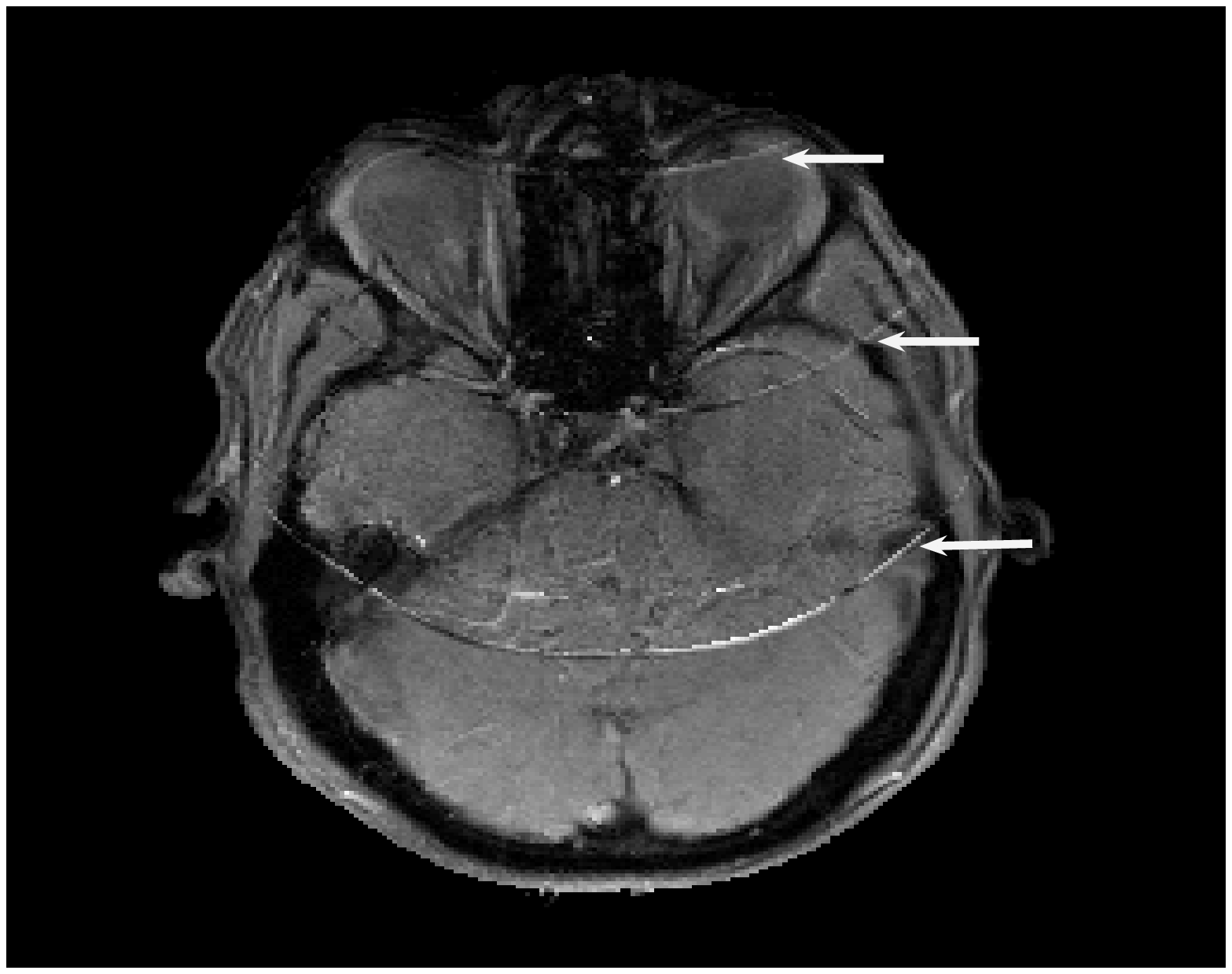} \\
\small SENSE & & \\
 & \includegraphics[width=4cm, height=4cm]{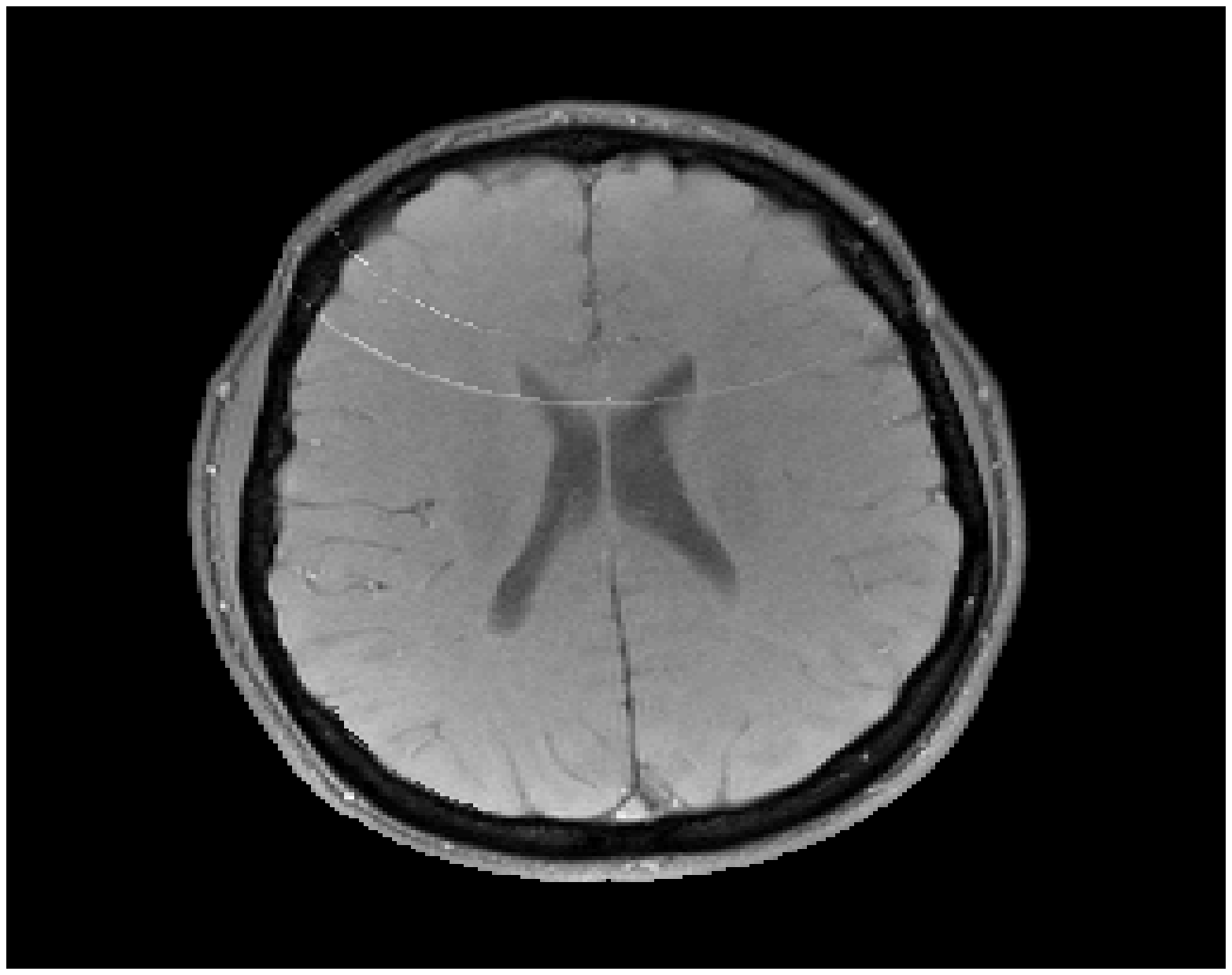} & \includegraphics[width=4cm, height=4cm]{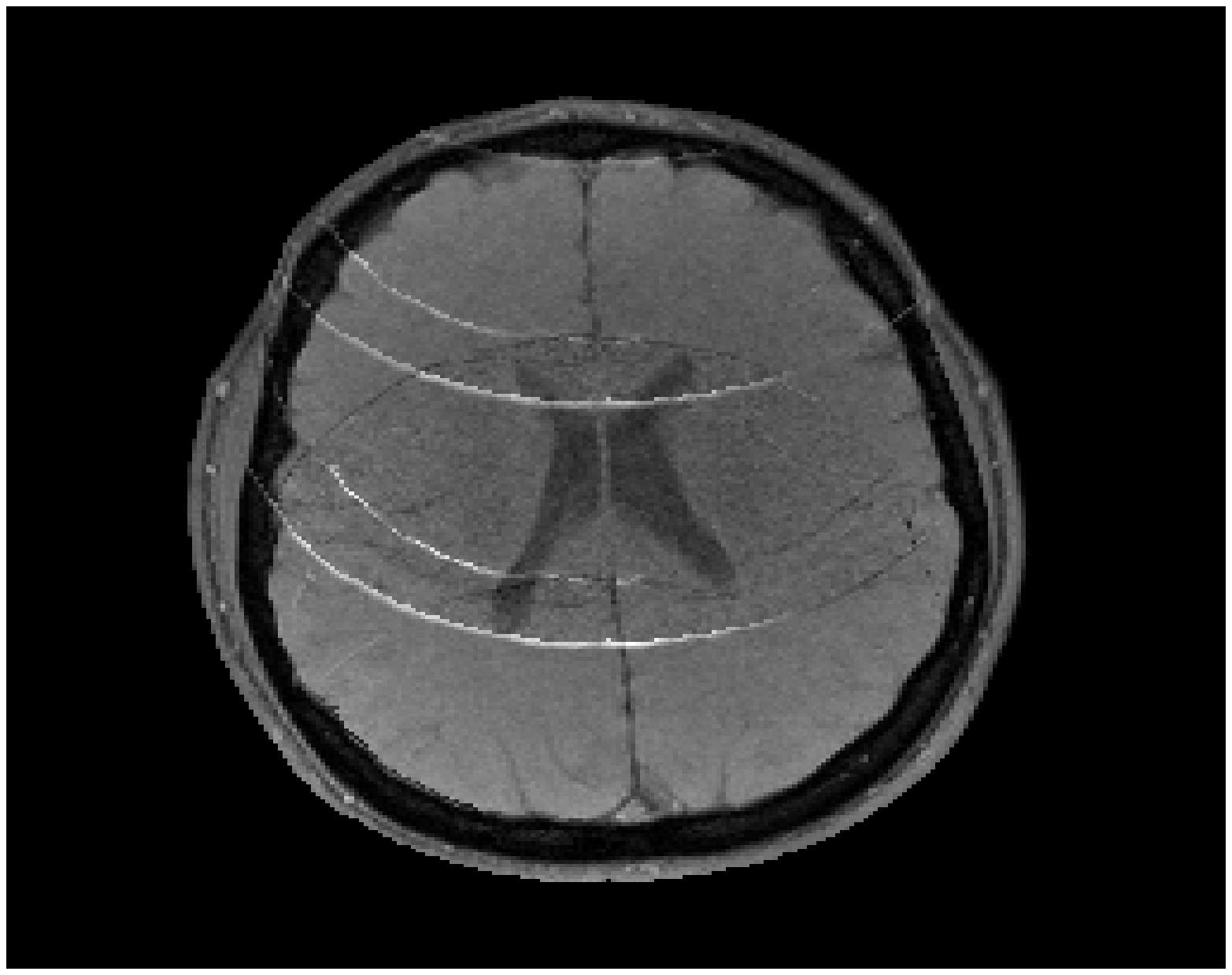} \\
\hline
\hline
 & \includegraphics[width=4cm, height=4cm]{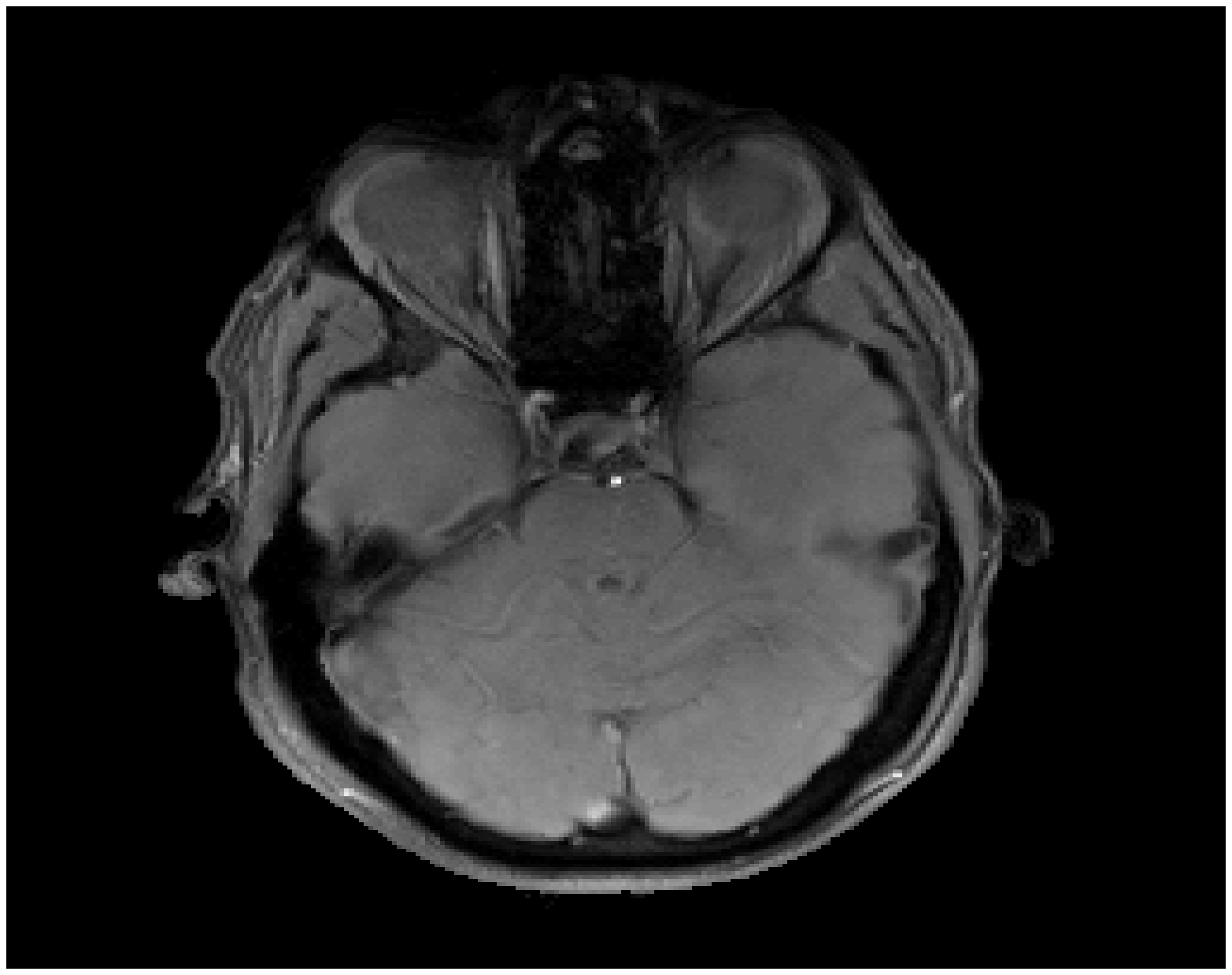} & \includegraphics[width=4cm, height=4cm]{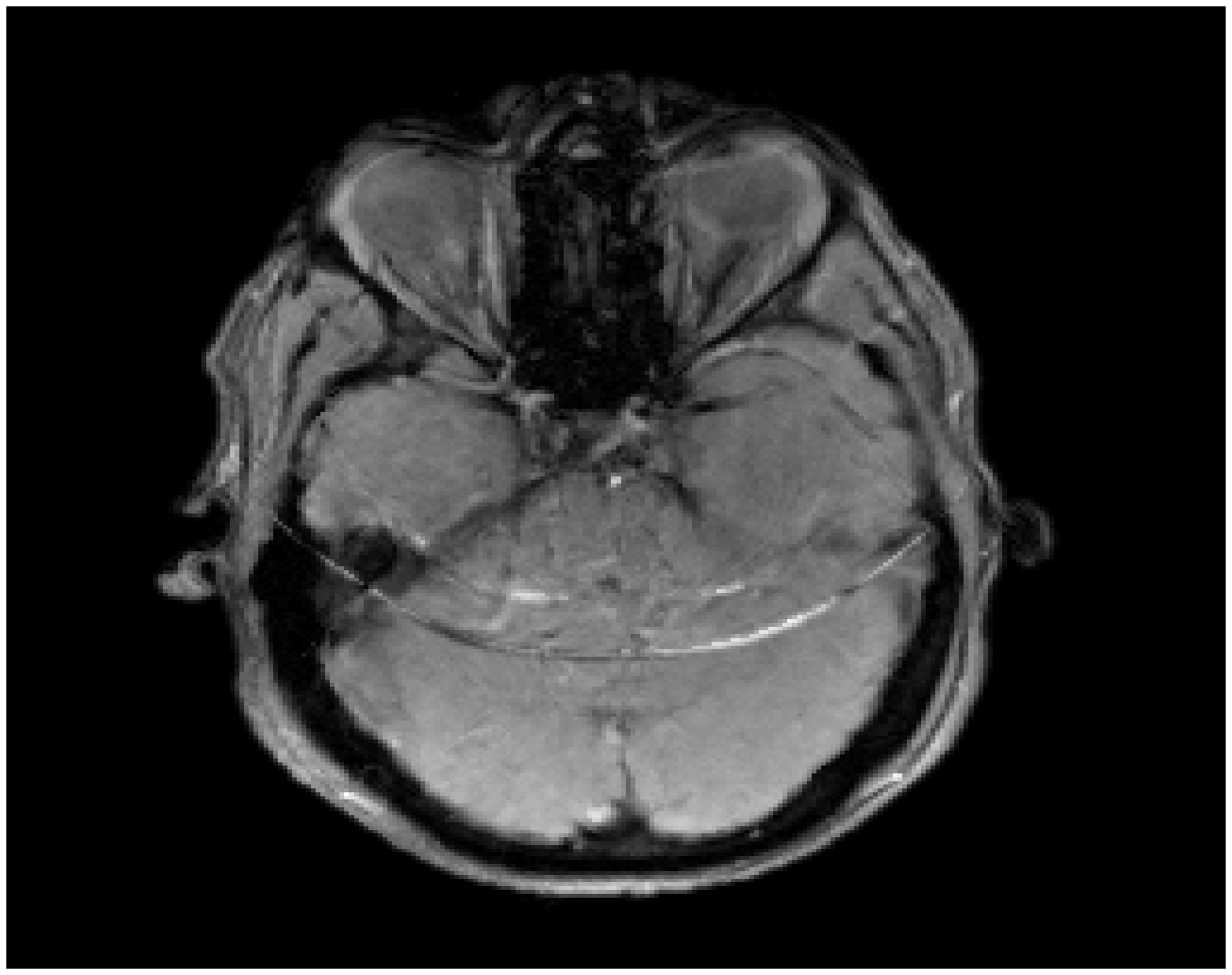} \\
\small Tikhonov & & \\
 & \includegraphics[width=4cm, height=4cm]{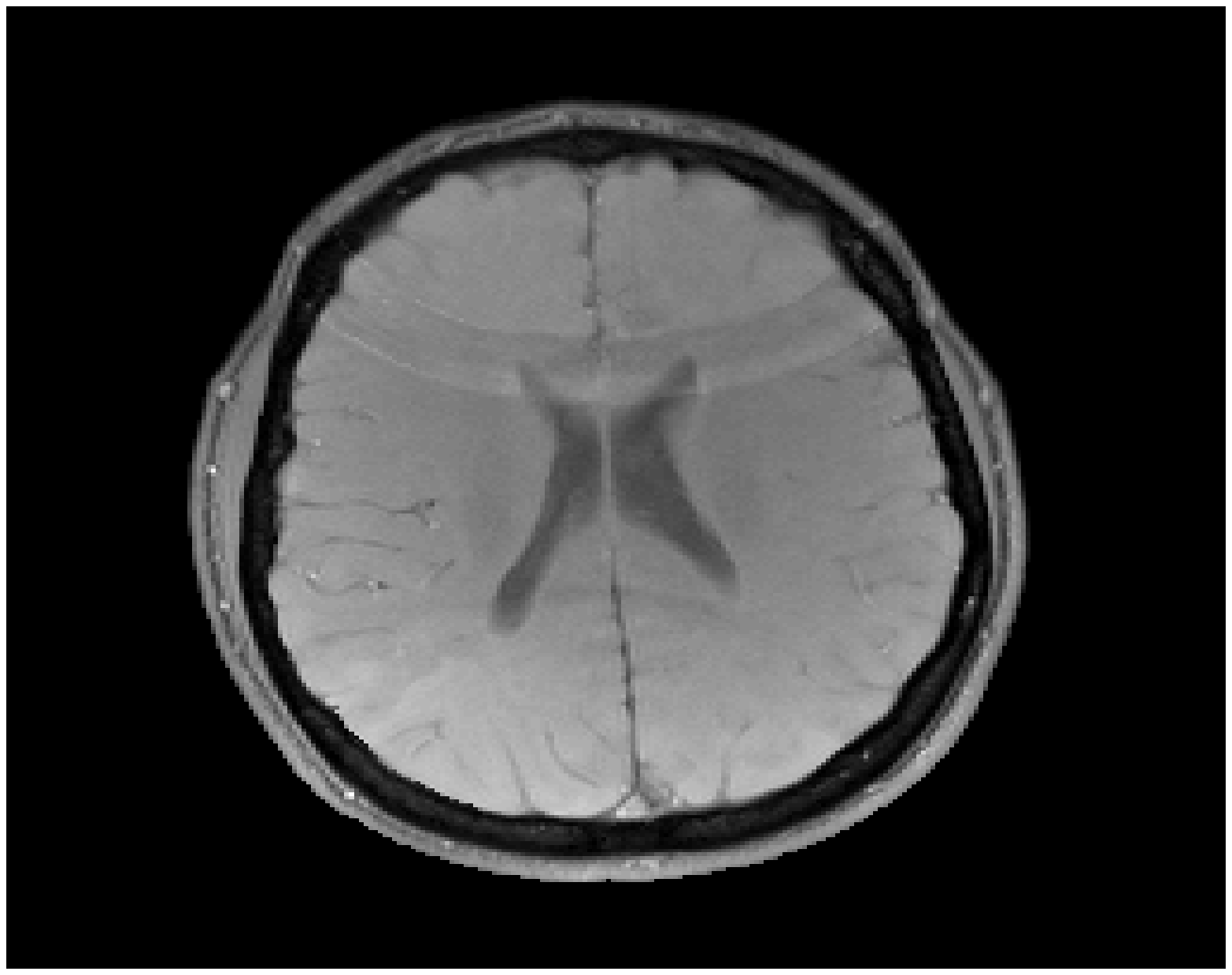} & \includegraphics[width=4cm, height=4cm]{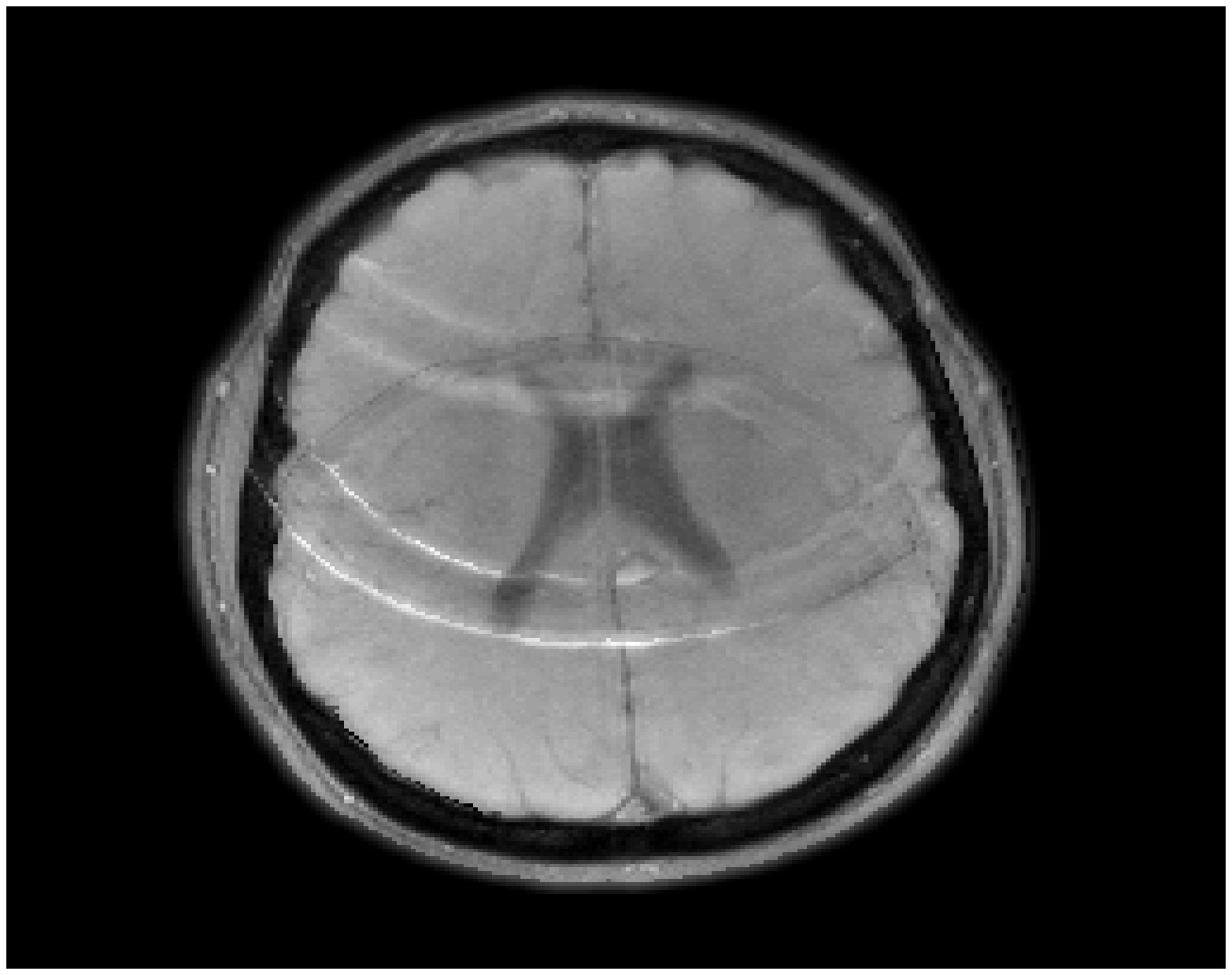} \\
\end{tabular}
\caption{Two reconstructed slices using SENSE and Tikhonov regularization for $R=2$ and $R=4$.}
\label{fig:sense_tikh}
\end{figure}

\begin{figure}
\centering
\vspace*{0.1cm}
\begin{tabular}{c c c}
& $R=2$&$R=4$\\
 & \includegraphics[width=4cm, height=4cm]{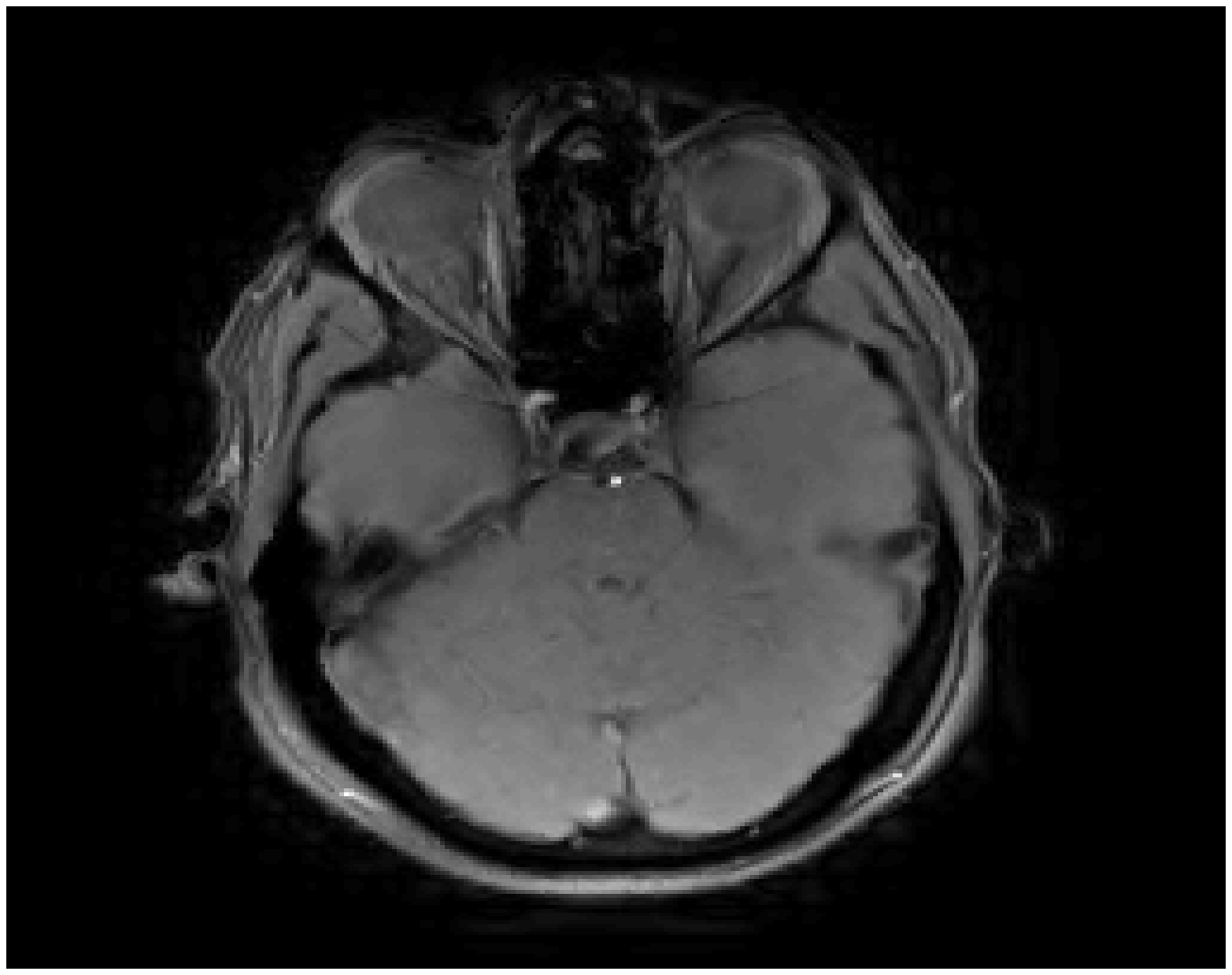} & \includegraphics[width=4cm, height=4cm]{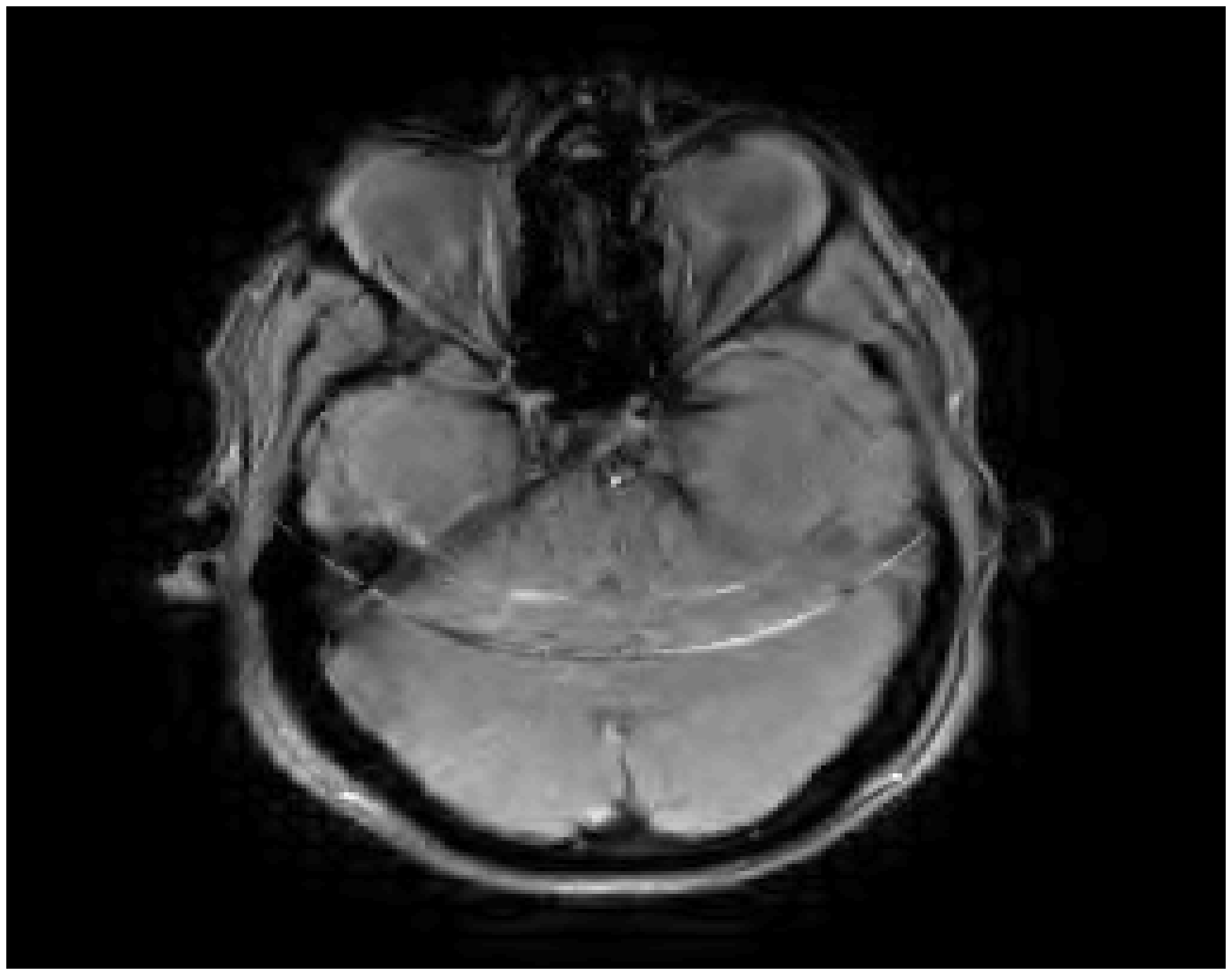} \\
\small Algorithm 1 & & \\
 & \includegraphics[width=4cm, height=4cm]{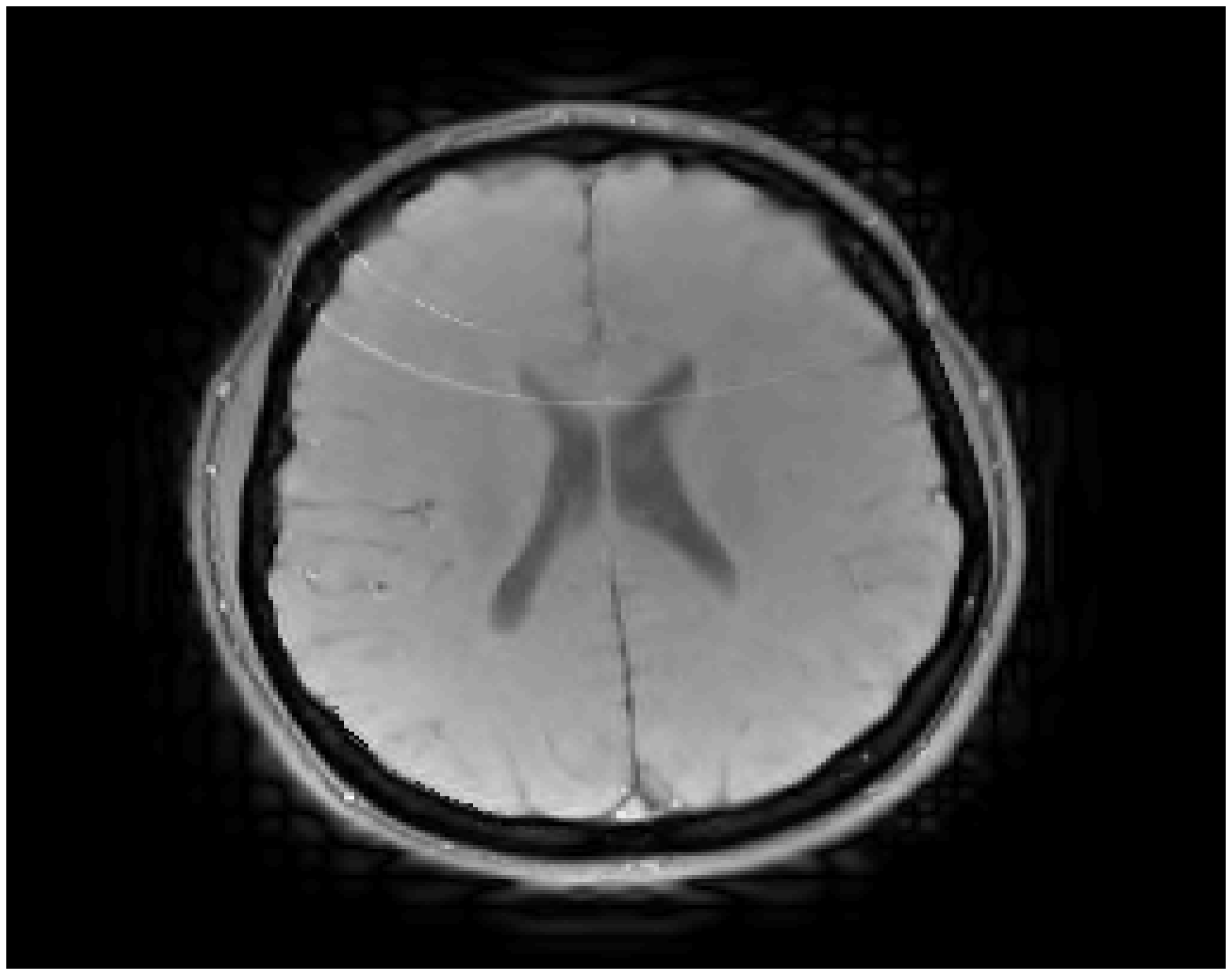} & \includegraphics[width=4cm, height=4cm]{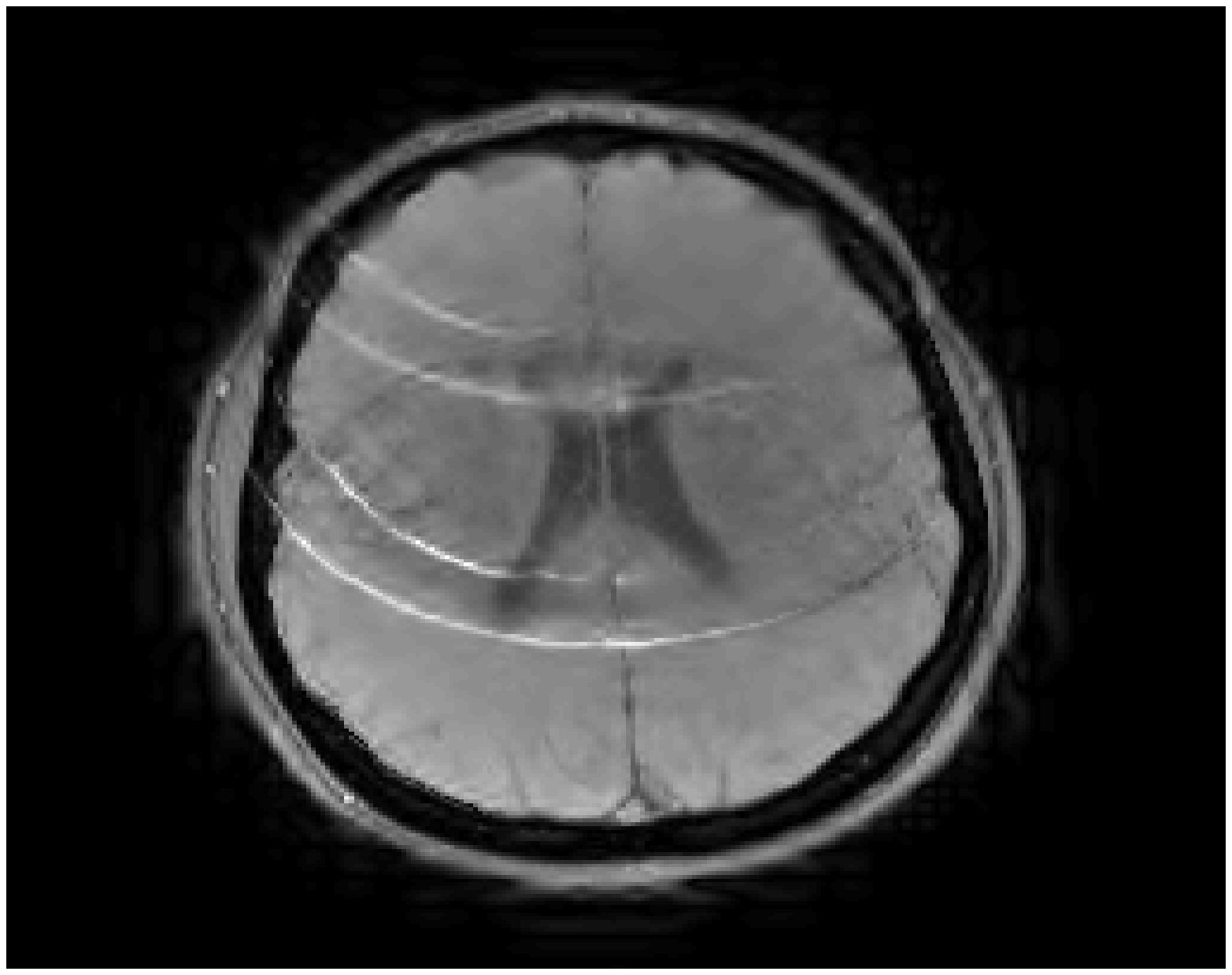} \\
\hline
\hline
 & \includegraphics[width=4cm, height=4cm]{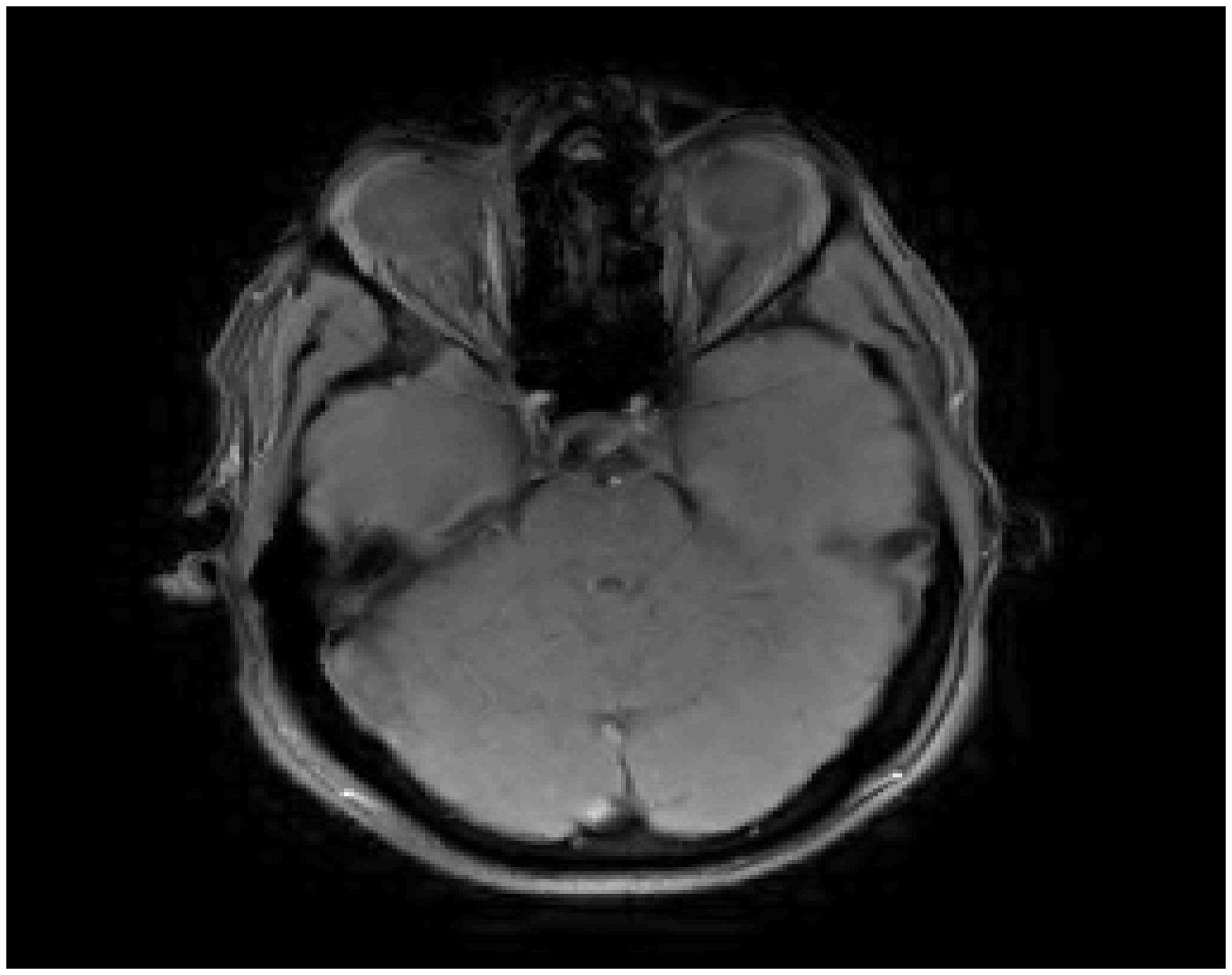} & \includegraphics[width=4cm, height=4cm]{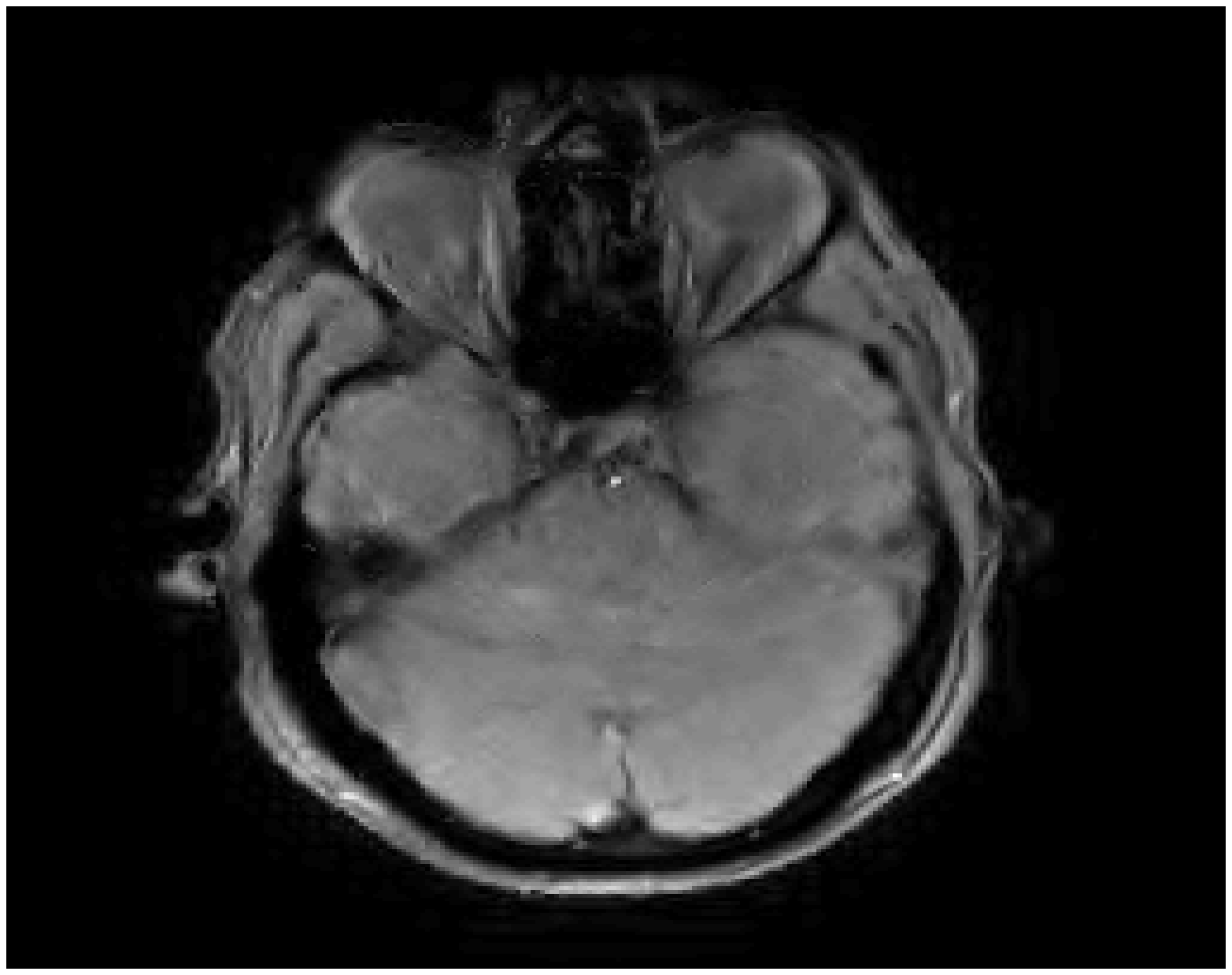} \\
\small Algorithm 2 & & \\
 & \includegraphics[width=4cm, height=4cm]{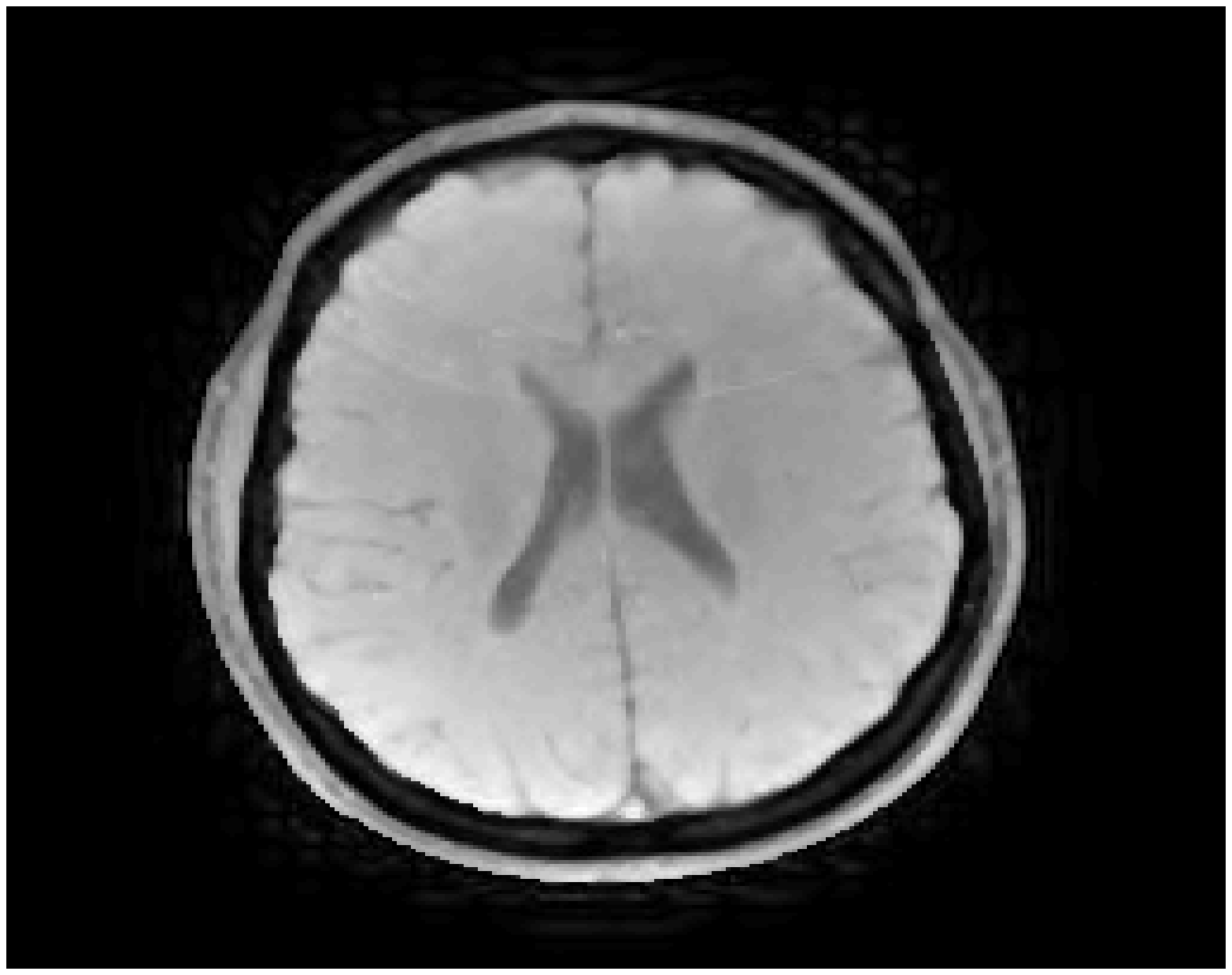} & \includegraphics[width=4cm, height=4cm]{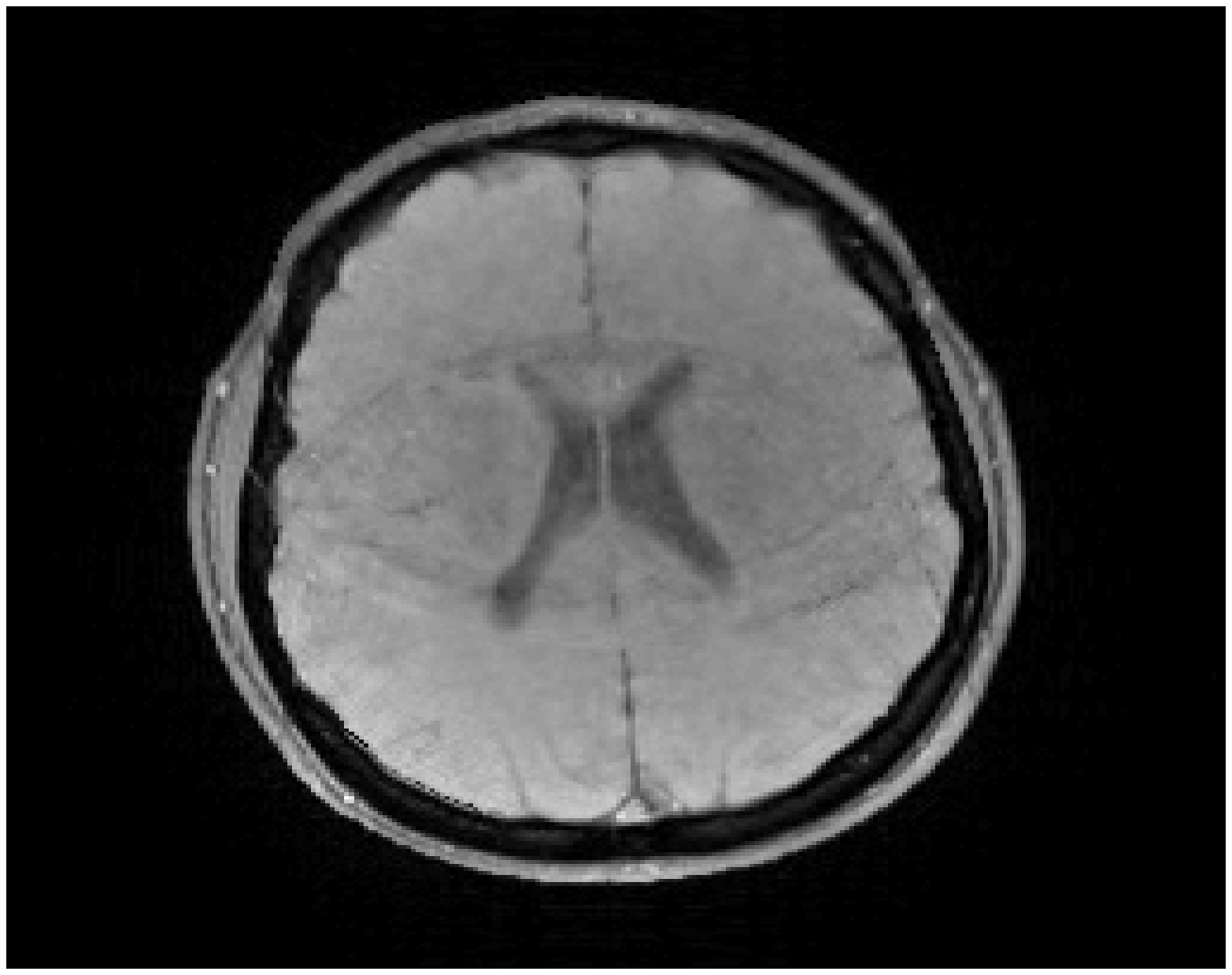} \\
\end{tabular}
\caption{Two reconstructed slices using Algorithm 1 and Algorithm 2 for $R=2$ and $R=4$.}
\label{fig:algos12}
\end{figure}

The algorithm has been stopped when the criterion no longer significantly varies, by choosing $\varepsilon =10^{-5}$ in Algorithm~\ref{algo:wt}. For different values of $\lambda$,  Fig.~\ref{fig:critere1} illustrates the evolution of the optimization criterion w.r.t the iteration number when reconstructing a 2D-slice. Comparison of the curves shows that $\lambda_n \equiv 1$ gives the fastest convergence.
One can also consider that, after about 20 iterations, which correspond to 160 seconds of execution time using Matlab 7.7 on an Intel Core 2 ($3$~GHz) architecture, the minimum was reached with $\lambda_n \equiv 1$, so providing the optimal MAP solution.\\
Note that accelerated algorithms such as TWIST or FISTA have been recently proposed in the literature \cite{Bioucas_07,Beck_09} for minimizing the same optimality criterion. However, for the considered application, it has been observed no significant improvement of the convergence profile by using these methods.\\
Fig.~\ref{fig:algos12} shows reconstructed full FOV anatomical images using the proposed approach (Algorithm~\ref{algo:wt}) for $R=2$ and $R=4$. To compare this reconstruction with the Tikhonov one, we can refer to Fig.~\ref{fig:sense_tikh}. Note that in Tikhonov regularization, the regularization image $\rho_\mathrm{r}$ was chosen as a mean image based on the basic-SENSE reconstruction, which contains the mean value of the signal of interest within the brain region. The regularization parameter $\kappa$ was manually fixed. A comparison with the basic-SENSE reconstruction can also be made from Fig.~\ref{fig:sense_tikh}.\\
It can be observed that the aliasing artifacts in the basic-SENSE reconstructed image are smoothed. They are completely removed in the areas where they were not actually very strong. When observing the reconstructed image using Tikhonov regularization, we see that it suffers from blurring effects (Fig.~\ref{fig:sense_tikh}). These drawbacks do no longer exist in the WT regularized image in Fig.~\ref{fig:algos12} where a good reconstruction is performed in the brain area.

If we compare the proposed approach with basic-SENSE and Tikhonov reconstructions from a quantitative viewpoint, clear improvements in the reconstruction process
are also obtained.
The quantitative measure we used is the Signal-to-Noise Ratio ($\mathrm{SNR}$).
Using the available reference image $\rho_{\mathrm{ref}}$ and for the reconstructed image $\hat{\rho}$ obtained with any reconstruction approach, it is possible to evaluate the signal-to-noise ratio: $\mathrm{SNR}=20\log_{10}\big(\|\rho_{\rm ref}\|/\|\rho_{\rm ref}-\hat{\rho}\|\big)$. Table \ref{tab:SNR} gives the $\mathrm{SNR}$ values corresponding to the basic-SENSE, Tikhonov regularization and the proposed WT regularization for different slices of the anatomical brain volume.
\begin{table}[!ht]
\centering
\caption{SNR evaluation for reconstructed images using different methods.}
\vspace*{0.1cm}
\begin{tabular}{|p{1.3 cm}|p{2 cm}|p{2.8 cm}|p{2.8 cm}|p{3cm}|}
\hline
&\multicolumn{4}{|c|}{$\mathrm{SNR}$ (dB)}\\
\cline{2-5}
& SENSE & Tikhonov regularization & WT regularization & Constrained WT regularization\\
\hline
Slice $1$ & 14.34 &14.48&14.67 & 15.53  \\
\hline
Slice $2$ & 11.53 &11.73&13.72 & 15.13  \\
\hline
Slice $3$ & 12.96 &13.37&14.02 & 14.12  \\
\hline
Slice $4$ & 9.22 &9.58&10.16 & 13.10  \\
\hline
Slice $5$ & 11.49 &11.88&12.06 & 12.25  \\
\hline
Slice $6$ & 9.67 &9.80&10.11 & 11.22\\
\hline
Slice $7$ & 11.04 &11.26&11.52 & 12.00  \\
\hline
Slice $8$ & 12.19 &12.60&12.92 & 13.60  \\
\hline
Slice 9  & 13.74 &13.28&14.29 & 15.66  \\
\hline
\end{tabular}
\label{tab:SNR}
\end{table}

 Note that reconstructed images here correspond to slices of the middle region of the brain volume for which full FOV images contain a large area of signal of interest w.r.t the image background. 

To improve reconstruction, the constrained algorithm presented in Section \ref{sec:constrained} was applied with the Symmlet 8 wavelet basis. The parameters $\lambda_{n}$ and $\gamma _{n}$ have been set to the same values as for the unconstrained case, while $\tau$ was fixed to 2 as it was practically observed that this value gives the best convergence rate for the underlying Douglas-Rachford iterations. In practice, the value of $M_n$ was defined as the minimal integer value such that $\vert \frac{\mathcal{\eta}^{(n,M_n-1)}-\mathcal{\eta}^{(n,M_n-2)}}{\mathcal{\eta}^{(n,M_n-2)}}\vert < 10^{-5}$, which results in about 4 iterations of the Douglas-Rachford algorithm.
A morphological gradient \cite{serra_89} was used to detect artifact regions on which we apply the additional convex constraint. The upper and lowed bounds which define the convex sets $C_\vect{r}$ (Eq. \eqref{eq:Cr}) were computed based on a morphological opening and closing applied to the basic-SENSE reconstructed image in order to discard very low and very high intensities. Reconstructed anatomical images using this constrained approach is displayed in Fig. \ref{fig:algos12}. 

Usual inspection of the results shows that the surviving artifacts in Fig.\ref{fig:sense_tikh} have now been removed. A comparison with basic-SENSE, Tikhonov and unconstrained WT regularizations can be  made from Figs. \ref{fig:sense_tikh} and \ref{fig:algos12}. Noticeable improvements in terms of $\mathrm{SNR}$ are  also obtained that can reach $3.6~\mathrm{dB}$ w.r.t basic-SENSE, $3.4~\mathrm{dB}$ w.r.t Tikhonov regularization and $2.94~\mathrm{dB}$ w.r.t unconstrained WT regularization (see Table \ref{tab:SNR}).
Finally, it is important to note that the performance of our constrained reconstruction also depends on the efficiency of the prior detection (based on morphological analysis) of the distorted area on which the convex constraints are applied.

\subsubsection{Choice of the wavelet basis}
\begin{figure}[!ht]
\centering
\begin{tabular}{c c}
\centering
\tiny
Symmlets 8 ($\mathrm{SNR}=15.66~\mathrm{dB}$)  & \tiny Daubechies 8 ($\mathrm{SNR}=15.60~\mathrm{dB}$)\\
\includegraphics[width=4cm, height=4cm]{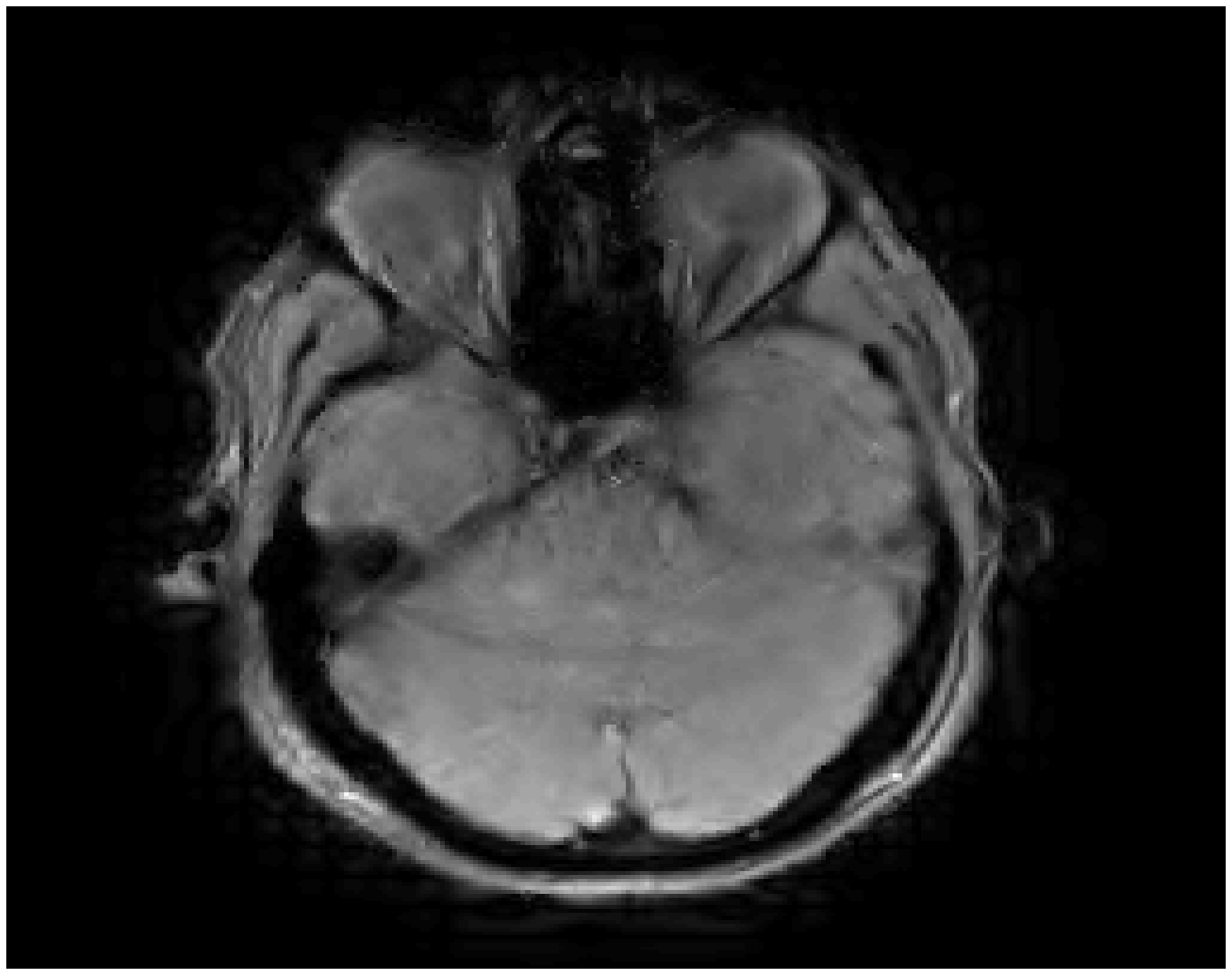} & 
\includegraphics[width=4cm, height=4cm]{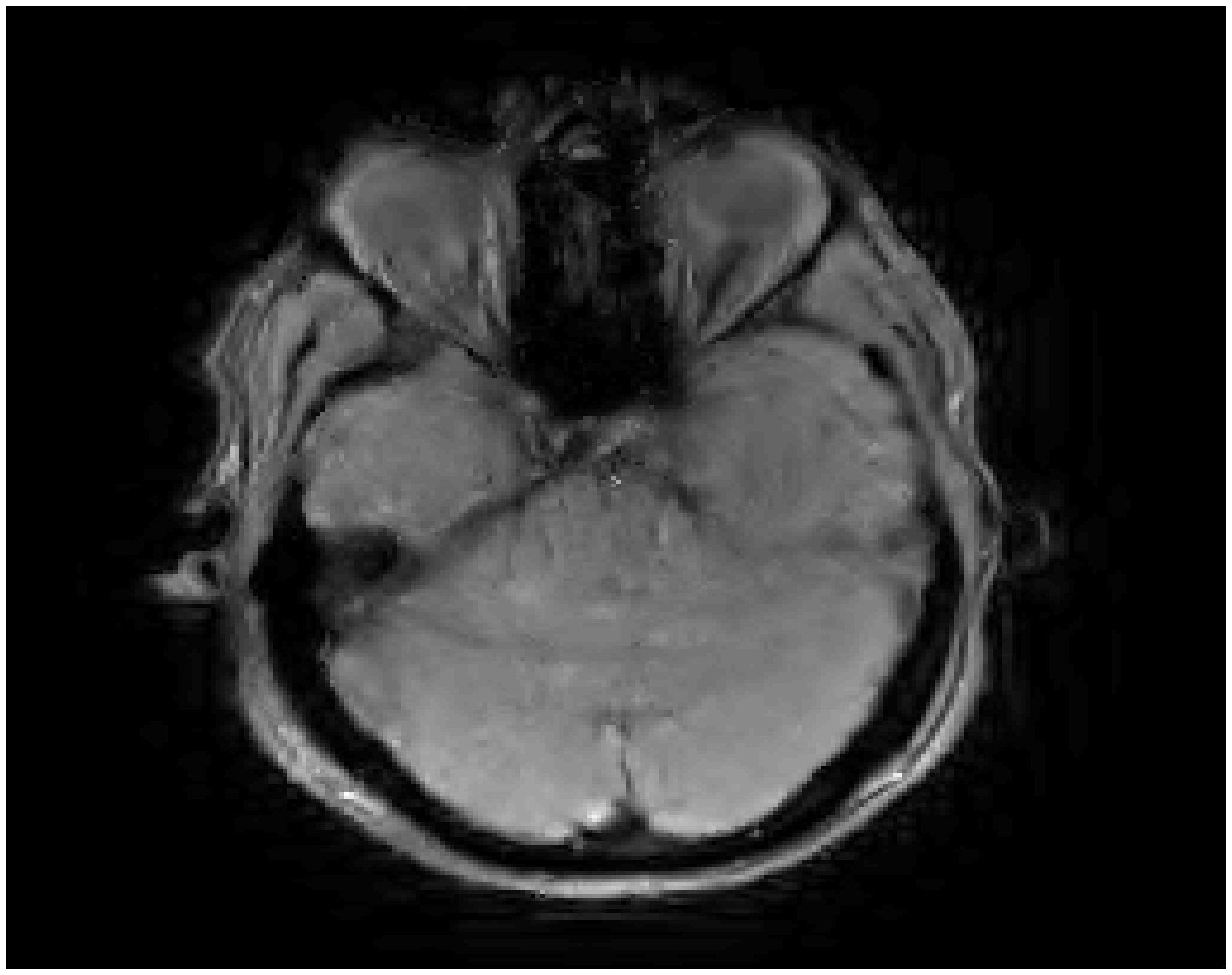} \\
\tiny
Haar ($\mathrm{SNR}=15.41~\mathrm{dB}$)  & \tiny Meyer 4 bands ($\mathrm{SNR}=12.56~\mathrm{dB}$)\\
\includegraphics[width=4cm, height=4cm]{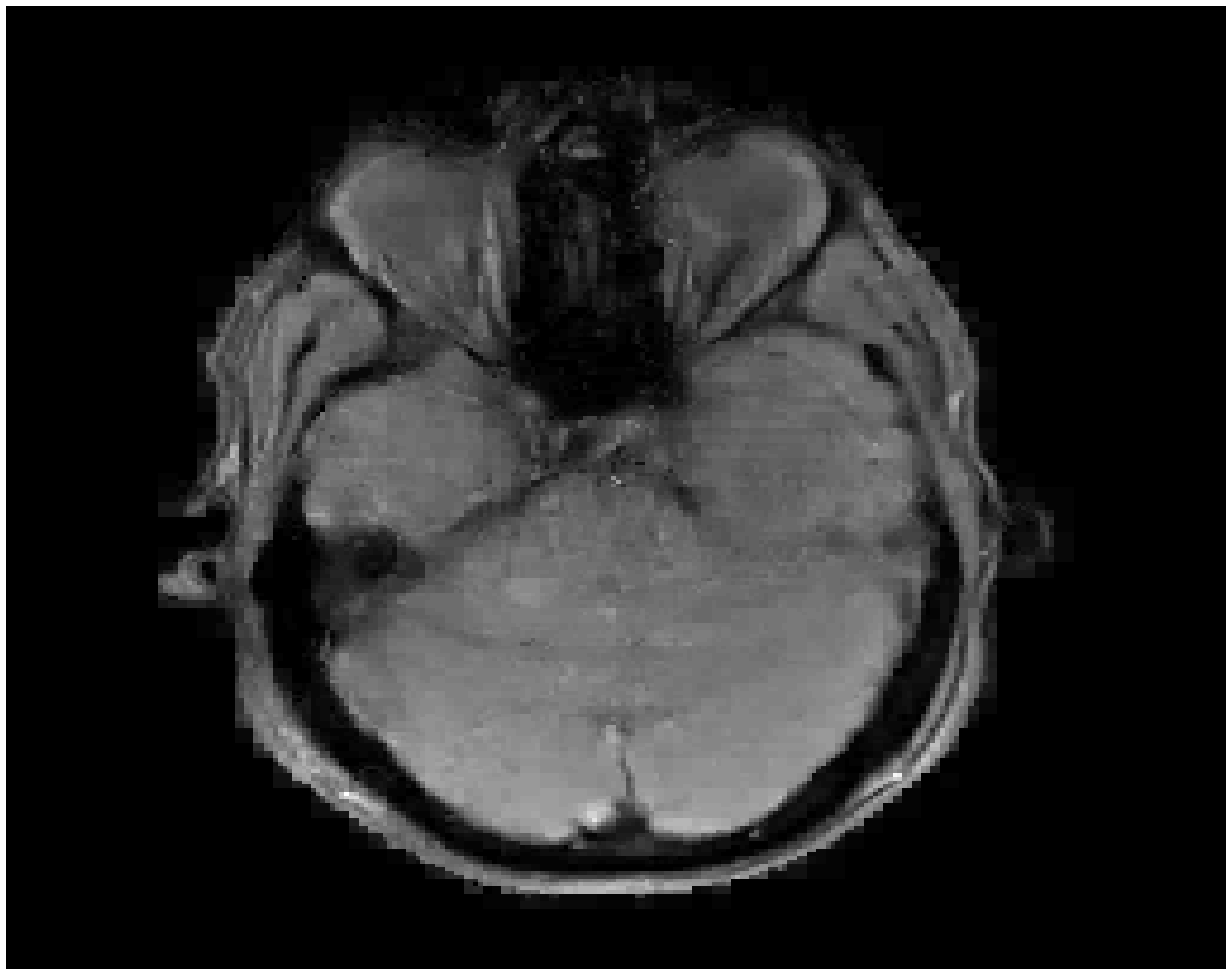} & 
\includegraphics[width=4cm, height=4cm]{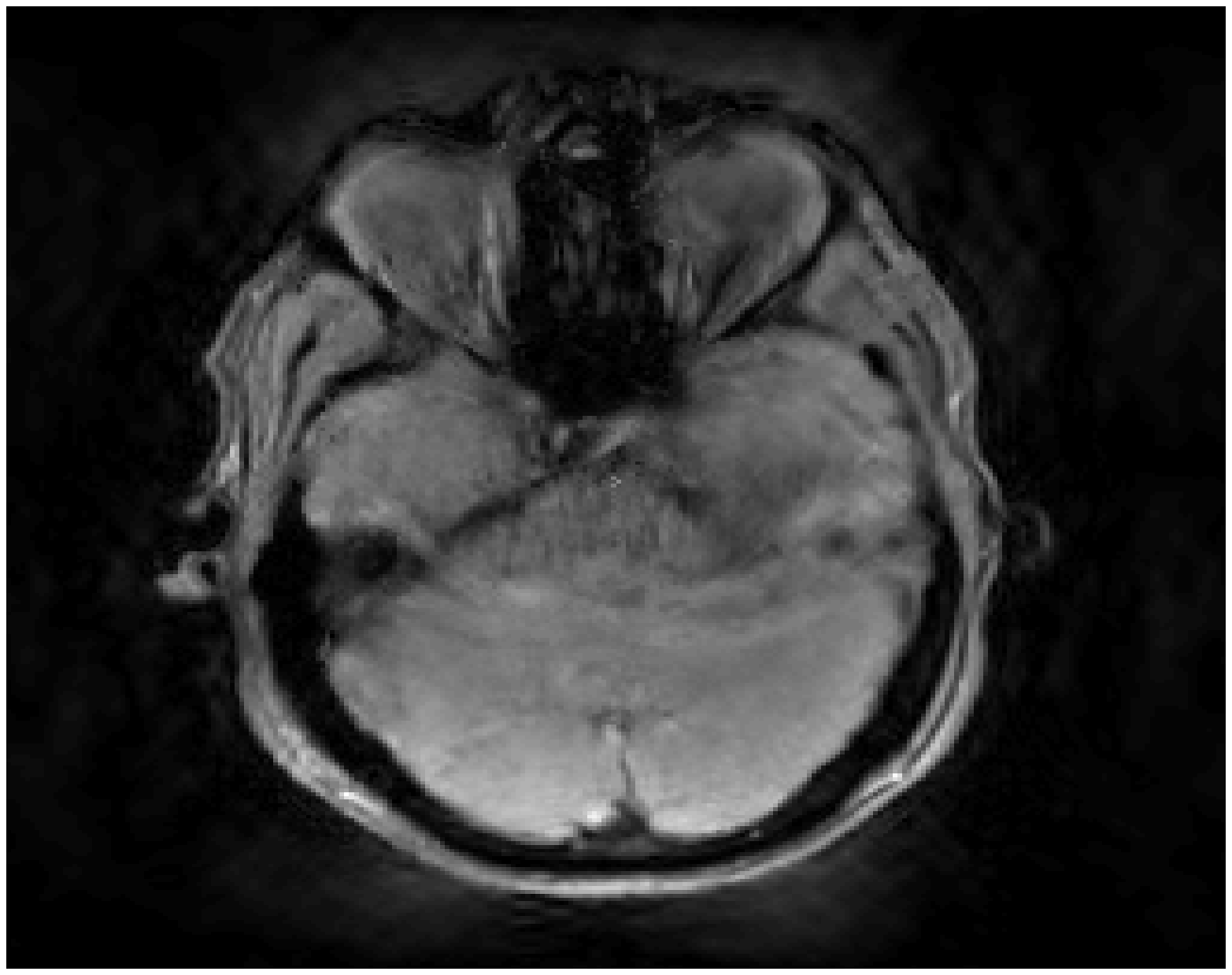} \\
\end{tabular}
\caption{Reconstructed images using different wavelet bases.}
\label{fig:diff_bases}
\end{figure}

\hfill \\
In this section, we study how the choice of the wavelet basis may influence the reconstruction performance. For comparison purposes, we present the results obtained with four different wavelet bases: dyadic Symmlet 8, dyadic Daubechies 8, dyadic Haar and Meyer with $M=4$ bands \cite{chaux_06}. In Fig. \ref{fig:diff_bases}, reconstructed images using the different wavelet bases with $j_{\mathrm{max}}=3$ are displayed.

Some boundary effects are introduced in the reconstructed images, but they are not very strong and do not affect the brain area of interest when using the Symmlets, Daubechies and Haar wavelet basis. However, these additional artifacts become very important when using the Meyer 4-band wavelet basis because of its large spatial support. Hence they drastically decrease the SNR of the reconstructed full FOV image. Note also that when using the Haar wavelet basis, we introduce some blocking effects caused by the Haar wavelet discontinuities that do not occur with the Symmlets and the Daubechies bases. Among the latters, Symmlets 8 gives slightly better regularization results than Daubechies 8.\\ 
More extensive comparisons using other bases were conducted, and we can conclude through this study that Symmlets 8 is the orthonormal wavelet basis among the tested ones which achieves the best WT regularization performance in terms of both $\mathrm{SNR}$ and visual quality. 
\subsubsection{Choice of the maximum resolution level} 
\begin{figure}[!ht]
\centering
\begin{tabular}{c c}
\centering
\tiny
(Image with $j_{\rm{max}}=2$)$-$(Image with $j_{\rm{max}}=1$)  & \tiny(Image with $j_{\rm{max}}=3$)$-$(Image with $j_{\rm{max}}=2$)\\
\includegraphics[width=4cm, height=4cm]{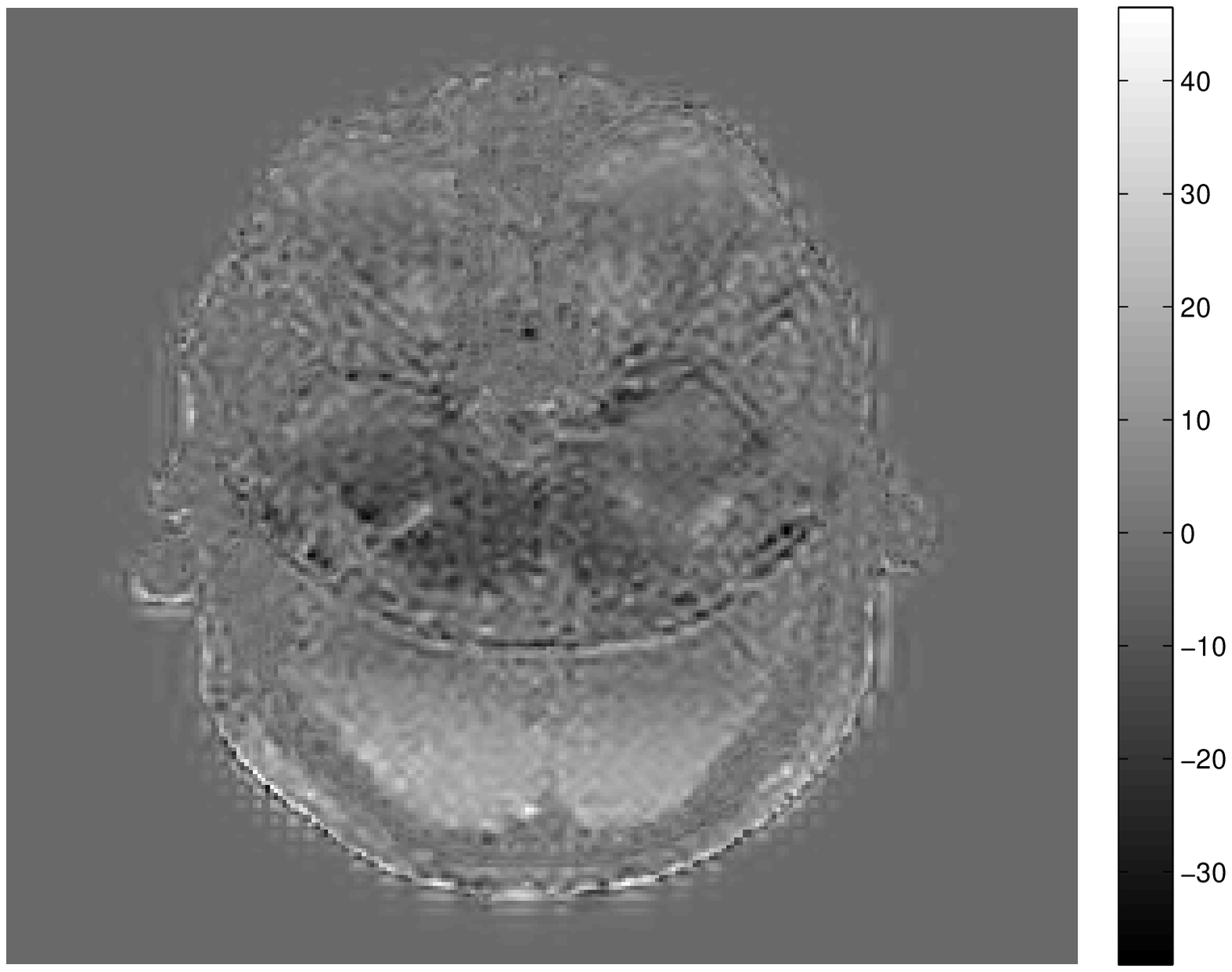} & 
\includegraphics[width=4cm, height=4cm]{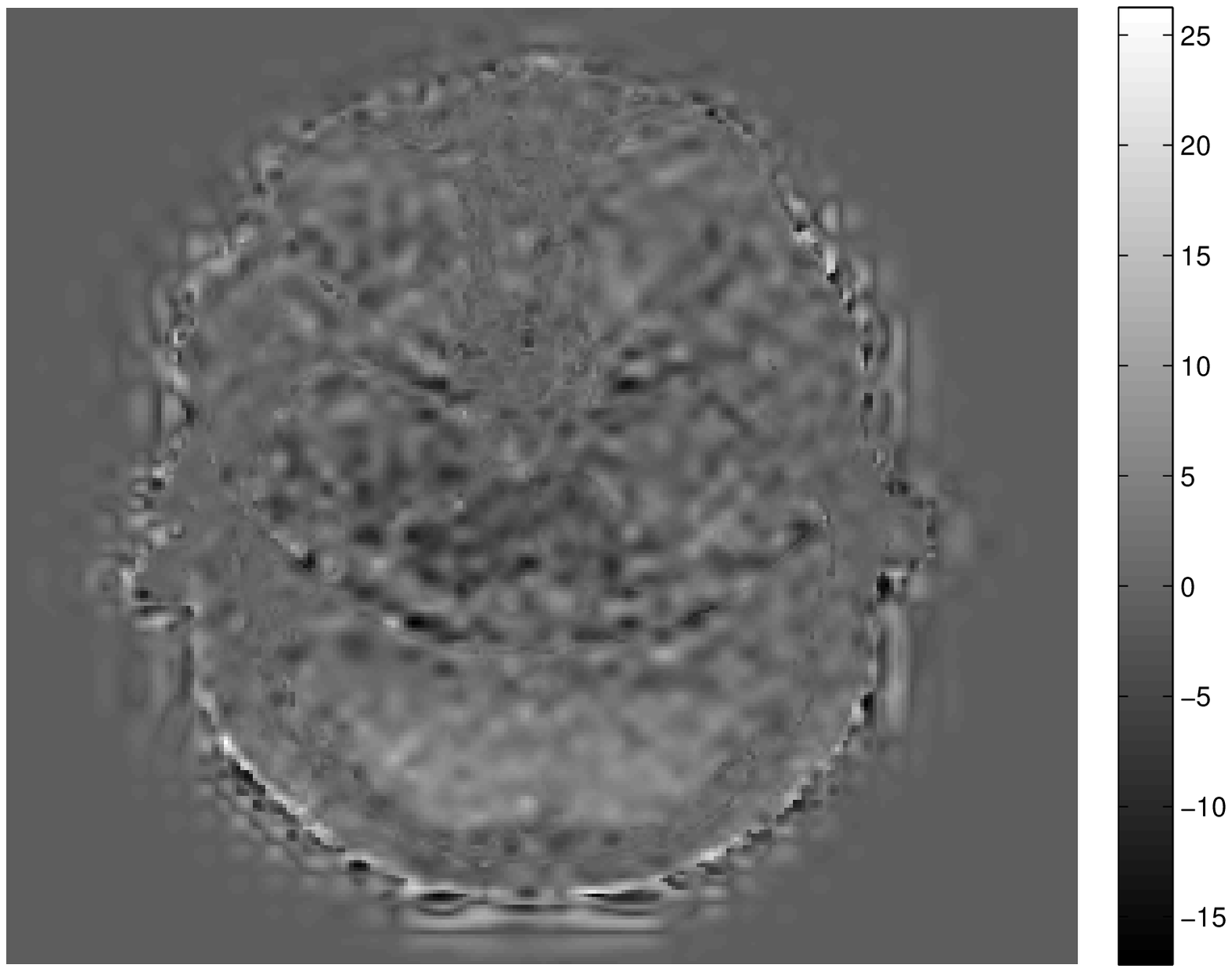} \\
\end{tabular}
\caption{Difference between reconstructed images using different resolution levels.}
\label{fig:diff}
\end{figure}

\hfill \\
In this section, we focus on the effect of the choice of the maximum resolution level $j_{\mathrm{max}}$ in terms of reconstruction quality.\\
The impact on reconstructed full FOV images can be emphasized through the difference between reconstructed images using 1, 2 and 3 resolution levels. Fig. \ref{fig:diff} illustrates the difference between reconstructed images using the first and the second (left), the second and the third (right) resolution levels. 

It appears that the difference between reconstructed images at different resolution levels is quite important since it can reach a value of 40 within an intensity range $[0,255]$. Moreover, this difference may be important in distorted areas. Hence, the higher the maximum resolution level, the better regularized the artifacts are.  However, only slight improvements are obtained beyond 3 resolution levels. \\
Note that by increasing the maximum resolution level, boundary effects are more visible, but they do not affect the brain.\\
A maximum resolution level $j_{\mathrm{max}}=3$ may therefore be a suitable choice for achieving an acceptable reconstructed image.

\subsection{Functional data}
\begin{figure}
\centering
\vspace*{0.1cm}
\begin{tabular}{c c c}
\small  SENSE & \includegraphics[width=3.5cm, height=3.5cm]{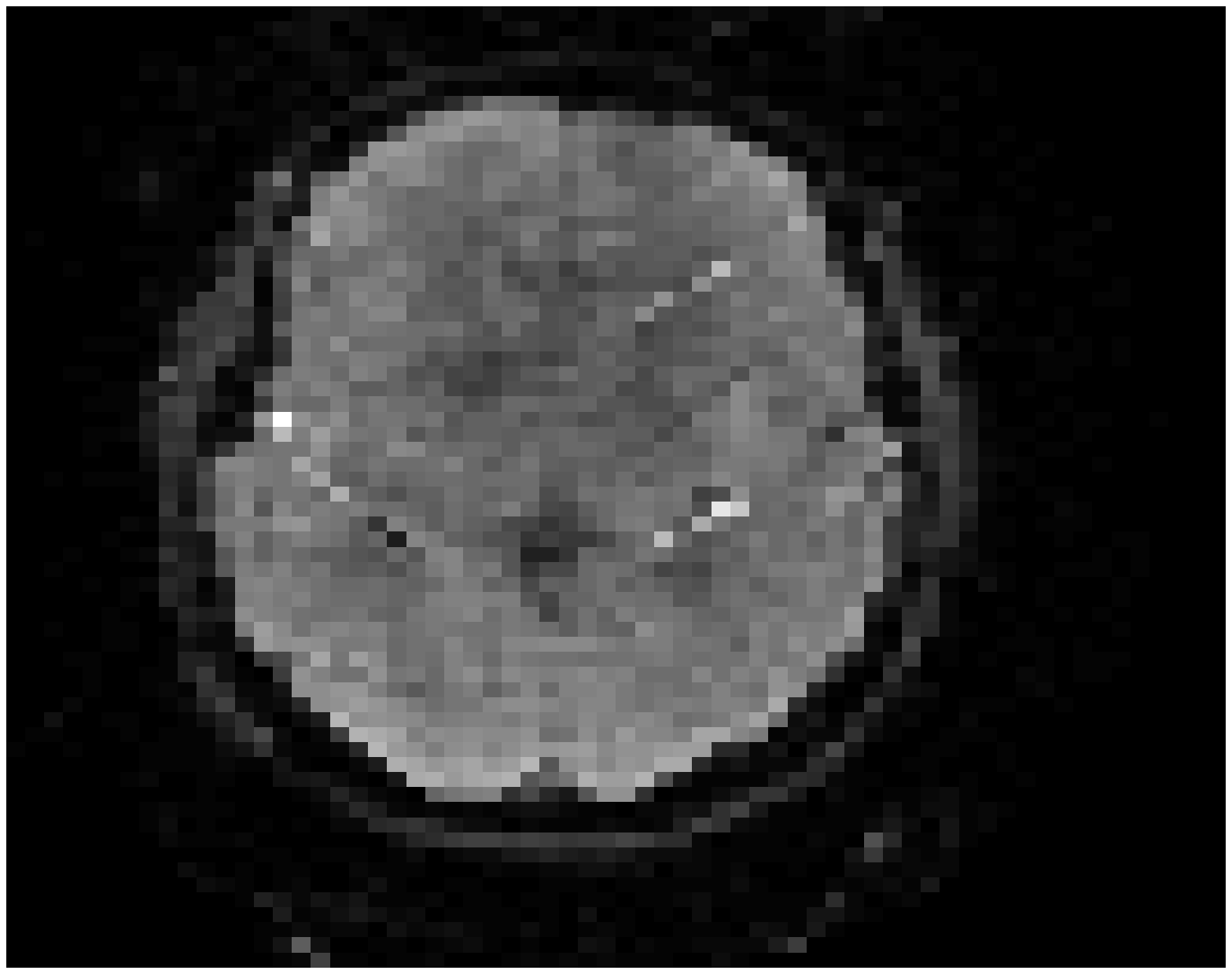} & \includegraphics[width=3.5cm, height=3.5cm]{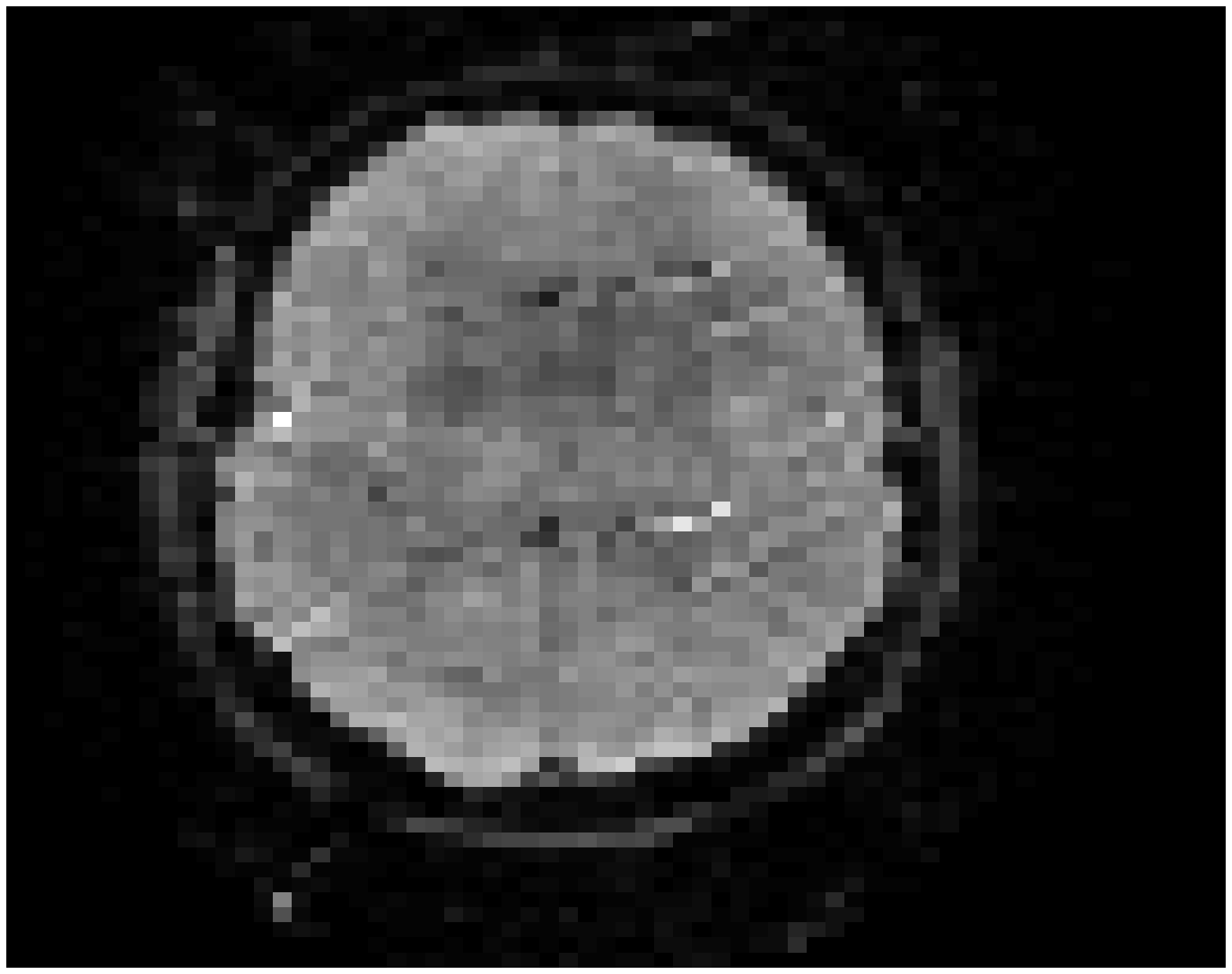} \\

\hline
\hline\\

\small Tikhonov  & \includegraphics[width=3.5cm, height=3.5cm]{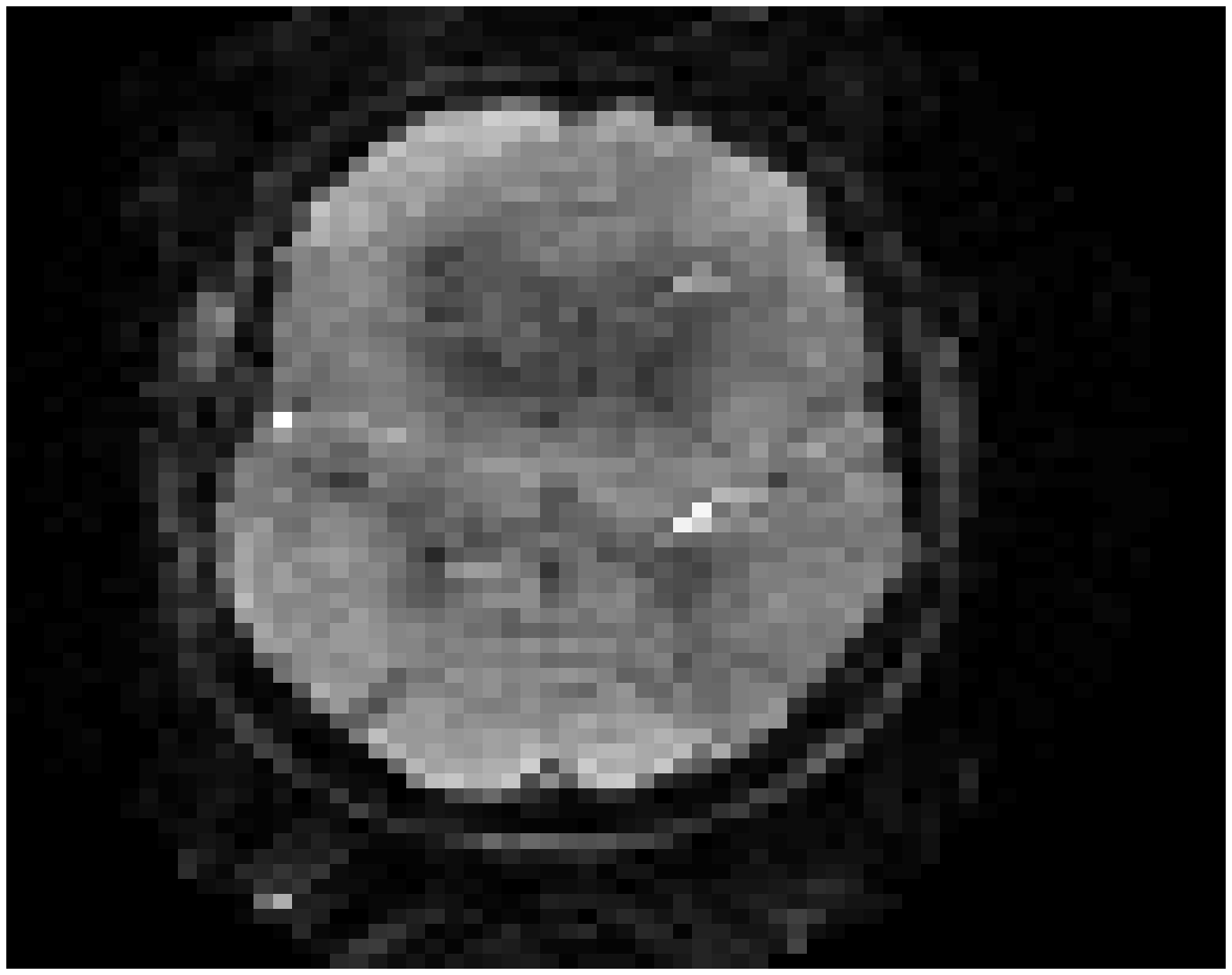} & \includegraphics[width=3.5cm, height=3.5cm]{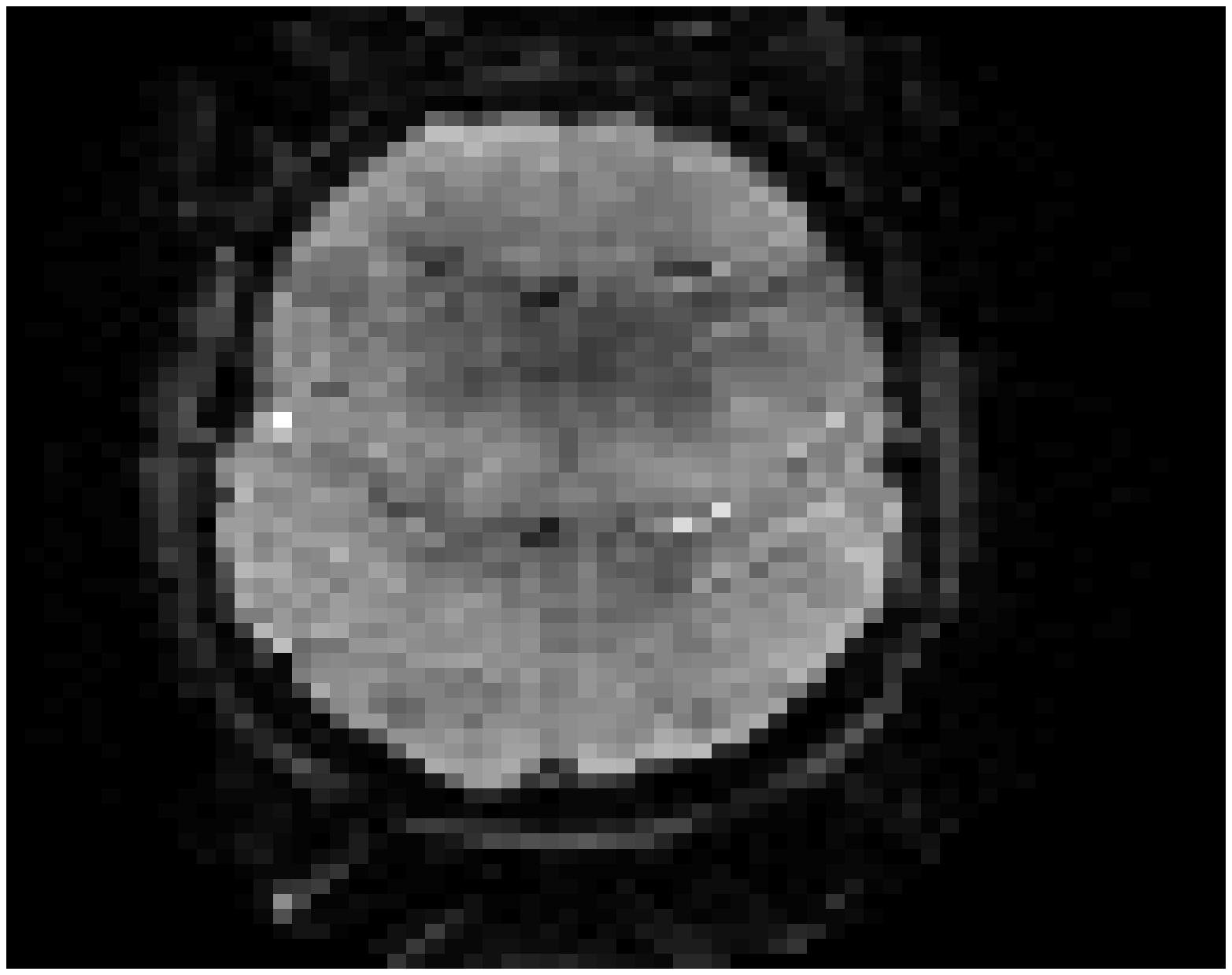} \\
\hline
\hline\\
\small Algorithm 1 & \includegraphics[width=3.5cm, height=3.5cm]{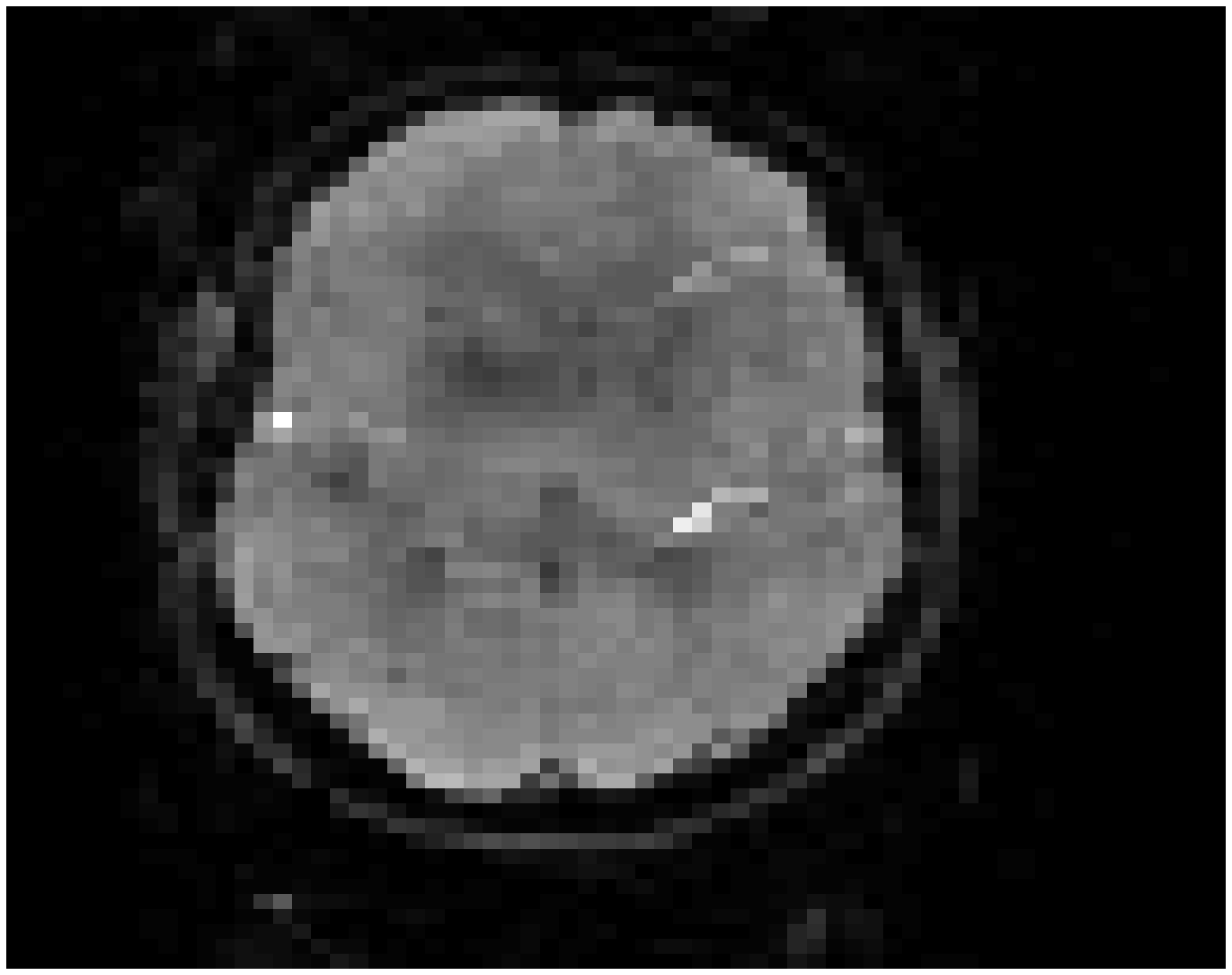} & \includegraphics[width=3.5cm, height=3.5cm]{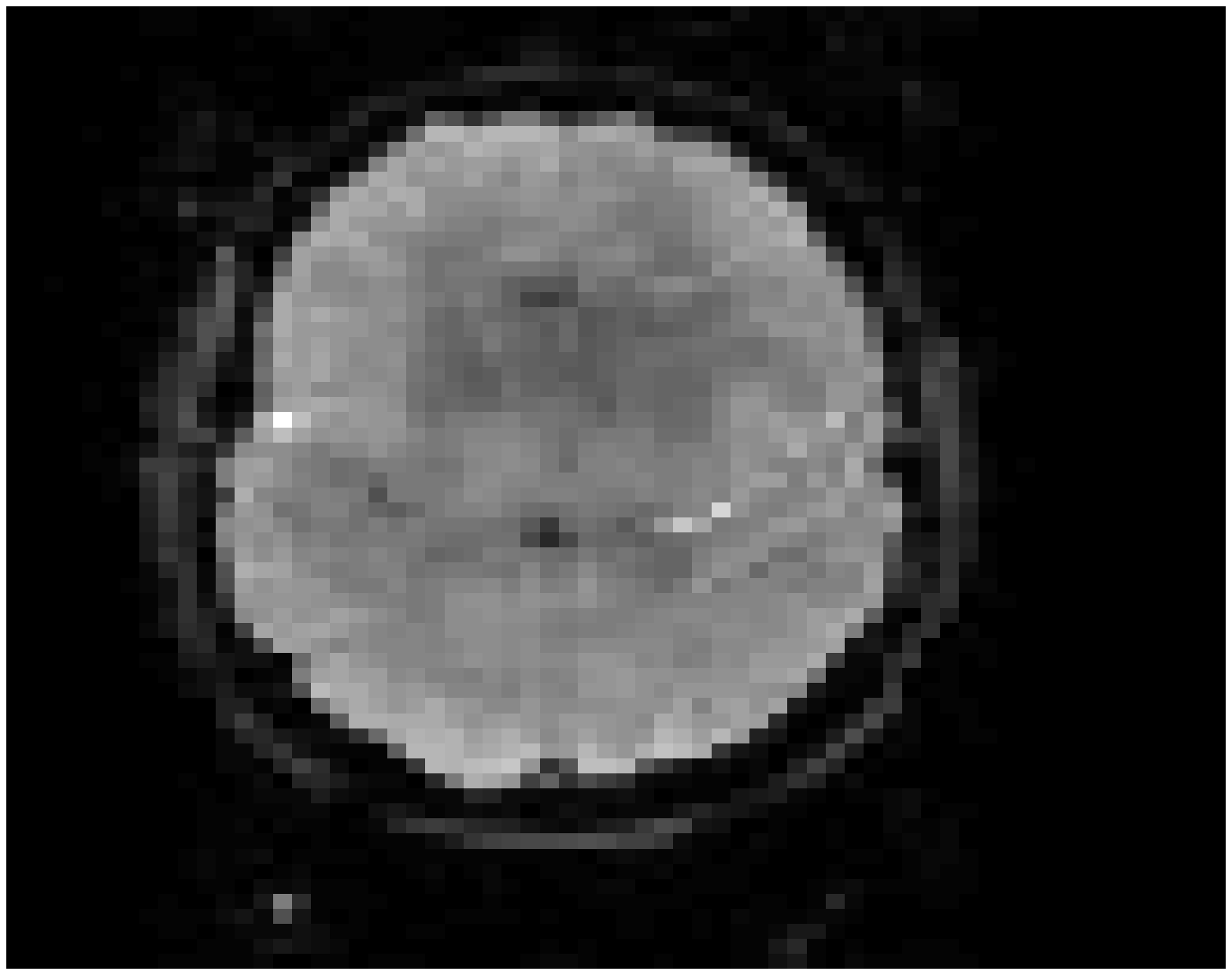} \\

\hline
\hline\\

\small Algorithm 2 & \includegraphics[width=3.5cm, height=3.5cm]{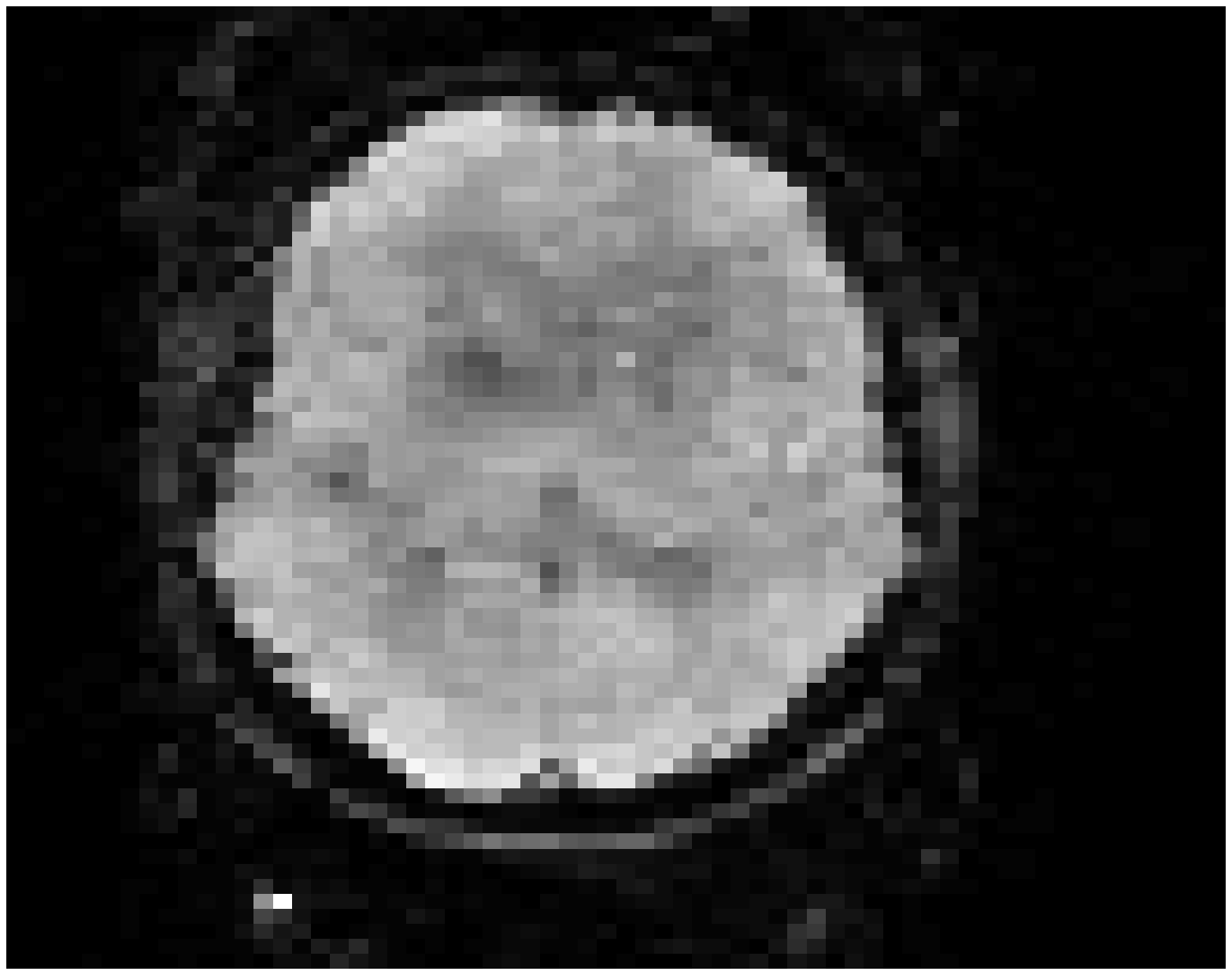} & \includegraphics[width=3.5cm, height=3.5cm]{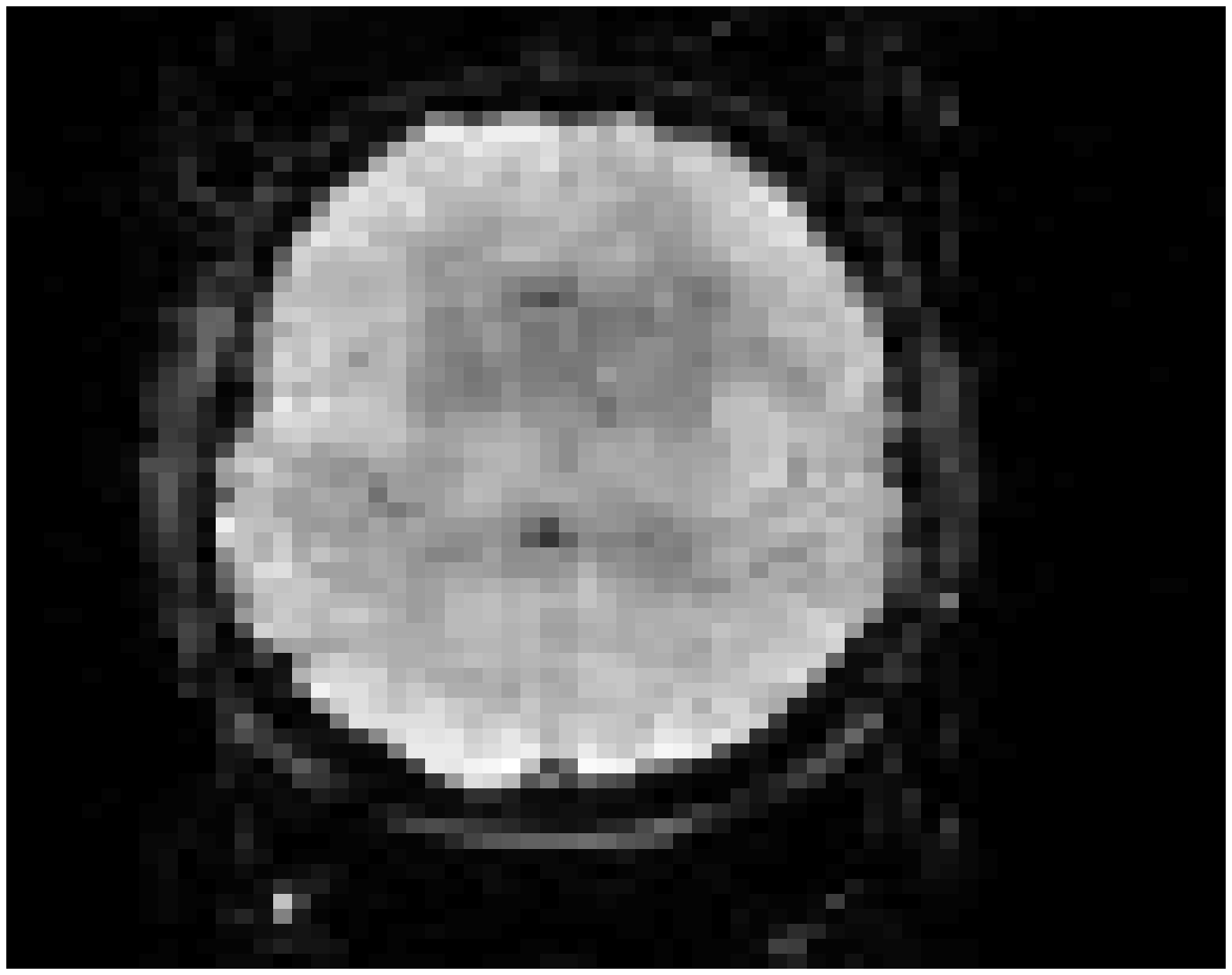} \\

\end{tabular}
\caption{Two reconstructed slices using SENSE, Tikhonov regularization, Algorithm 1 and Algorithm 2 for $R=4$.}
\label{fig:sense_tikh_f}
\end{figure}

This experiment on functional EPI data of size $64\times64$ was conducted using the same wavelet basis and priors. Algorithm parameters (i.e. relaxation and step-size parameters) have been adjusted according to the same rules as for anatomical data. The relaxation parameter was set to 1 and the step-size parameter was set to $1.99/\theta = 20.63$ after computation of the $\theta$ constant. Note that these EPI images have been acquired with no experimental paradigm at resting state (eyes closed and subject lying on the MRI bed).
Fig. \ref{fig:sense_tikh_f} illustrates reconstructed full FOV slices using SENSE, Tikhonov, Algorithm 1 and Algorithm 2.

It can be shown that many defective pixels were corrected when using the proposed WT regularization in Algorithm 1 without introducing additional artifacts, in contrast with Tikhonov regularization.
When analyzing reconstructed images using Algorithm 2, it can be noticed that pixels with very high intensity were smoothed, when compared with Tikhonov and basic-SENSE reconstructed images. Residual defective pixels belonging to distorted areas in SENSE reconstructed image have also been removed due to the convex constraint introduced in these areas. Note that the same approach as that applied in anatomical experiments has been adopted to detect artifact areas and to compute the upper and lower bounds defining the convex set $C_{\vect{r}}$. 
\subsubsection{Analysis of resting state BOLD signal}
\begin{center}
\begin{figure}[!ht]
\centering
\includegraphics[width=13cm, height=4cm]{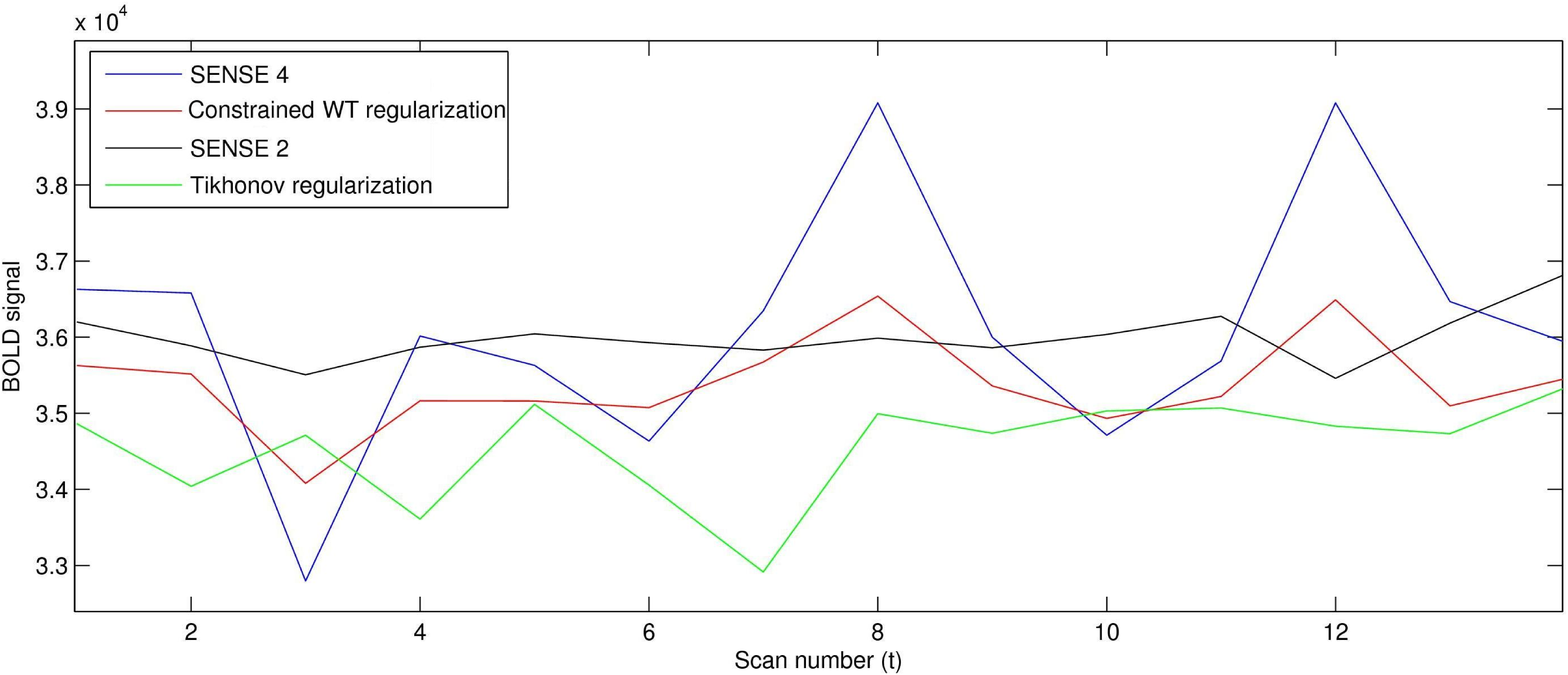}
\caption{Time evolution of the BOLD signal in the occipital cortex.}
\label{fig:stat}
\end{figure}
\end{center}

\hfill \\
To validate our reconstruction approach, an analysis was carried out on the occipital cortex area, where a statistical study is possible even if we have resting-state images with no experimental paradigm. Fig.~\ref{fig:stat} shows the variation of the BOLD signal w.r.t. time for images reconstructed with $R=4$ using SENSE, Tikhonov and the proposed approach. The variation of the BOLD signal of images reconstructed using SENSE with $R=2$ is also given and considered as a reference signal since with $R=2$, SENSE provides a quite reasonable reconstruction quality.

The four plots indicate that, by using our WT approach (constrained wavelet-based regularization), the BOLD signal is smoother than when using SENSE or Tikhonov regularization for EPI images acquired with $R=4$. In fact, the occipital cortex is a brain region which is always activated even if there is no experimental paradigm, which justifies the smoothness of the BOLD signal in this region. 

\section{Discussion and Conclusions}\label{sec:concl}
In this paper, we proposed and tested a novel approach for SENSE reconstruction based on a regularization in a 2D orthonormal wavelet basis. This method reduces aliasing artifacts related to the noise in the acquired data, which become critical when using high reduction factors with low magnetic field intensity. Our results on 1.5 Tesla data show improvements in reconstructed images compared with basic SENSE and other standard regularization methods for both anatomical and functional images. Different choices for the wavelet basis and the maximum resolution level have been discussed to evaluate the impact of these configurations on our algorithm performance. Another potential advantage of the proposed method is that it can be implemented with parallel CPU architectures as we can deal with each wavelet coefficient subband separately in some steps of the algorithm. Parallelism can also be exploited to reconstruct different slices or, in fMRI, different volumes independently. Since our approach is more time-consuming than SENSE when reconstructing a 2D-slice, this parallelism let us catch up computational time when reconstructing whole brain volumes or temporal series of brain volumes.\\
A constrained version of our algorithm was also presented in order to achieve more accurate reconstructed images. Our results show improvements both in visual image quality and in quantitative measures. \\
Concerning our future work, a more sophisticated detection of artifact regions could improve the obtained results. The incorporation of additional convex constraints on artifact regions could also be of great interest. An extension to 3D wavelet decompositions may also be envisaged to reconstruct whole brain volumes at once and probably exploit between-slices spatial dependencies. Experimentations of our approach on 3 Tesla data and functional data with experimental paradigm may also be important.

\appendix 
\begin{center}\label{append:a1}
    APPENDIX A\\ 
Proof of Proposition \ref{prop:cvl}
\end{center}

From \eqref{eq:J20}-\eqref{eq:J22}, it can be seen that $\mathcal{J}_P$
is a convex function such that
\[
\forall \zeta \in \mathbb{C}^K,\quad
\mathcal{J}_P(\zeta) \ge \frac{\vartheta_1}{2} \|\zeta\|^2-\vartheta_0
\]
where 
\begin{align}
\vartheta_0 &= \frac{K_{j_{\mathrm{max}}}}{2} \Big(\frac{(\mu^{\mathrm{Re}})^2}{\sigma_{\mathrm{Re}}^{2}}+\frac{(\mu^{\mathrm{Im}})^2}{\sigma_{\mathrm{Im}}^{2}}\Big) \nonumber \\ 
\vartheta_1 &= \min\{(2\sigma_{\mathrm{Re}}^{2})^{-1},(2\sigma_{\mathrm{Im}}^{2})^{-1},
(\beta_{j,o}^{\mathrm{Re}})_{o \in \mathbb{O},1\le j \le j_{\mathrm{max}}},
(\beta_{j,o}^{\mathrm{Im}})_{o \in \mathbb{O},1\le j \le j_{\mathrm{max}}}\}.\nonumber
\end{align}
This means that $\mathcal{J}_P$ is a strongly convex function with modulus
$\vartheta_L$ \footnote{A function $f\colon \chi \to ]-\infty,+\infty]$,
where $\chi$ is a Hilbert space, is said \textit{strongly convex} on $\chi$ with modulus $\vartheta_1 > 0$ if there exists some $g \in \Gamma_0(\chi)$ such that $f=g+\frac{\vartheta_1\|\cdot\|^2}{2}$.}. Since $\mathcal{J}_L$ is a finite convex function,
$\mathcal{J}_{\rm WT}$ also is strongly convex. It is thus strictly convex
and coercive (i.e. $\lim_{\|\zeta \| \to \pinf} \mathcal{J}_{\rm WT}(\zeta) = \pinf$)
and, from standard result in convex analysis \cite{ekeland_99}, it can be deduced
that $\mathcal{J}_{\rm WT}$ has a unique minimizer $\widehat{\zeta}$.

In addition, $\mathcal{J}_L$ is a differentiable function and
we have
\begin{equation}
\forall \zeta \in \mathbb{C}^K,\qquad
\nabla \mathcal{J}_{LT}(\zeta) 
= \frac{\partial \mathcal{J}_{LT}(\zeta)}{\partial \mathrm{Re}(\zeta)}
+ \imath \frac{\partial \mathcal{J}_{LT}(\zeta)}{\partial \mathrm{Im}(\zeta)}
= T\, \nabla \mathcal{J}_L(T^*\zeta),
\end{equation}
where $\mathcal{J}_{LT}(\zeta) = \mathcal{J}_{L}(T^*\zeta)$.
Set $\rho= T^*\zeta$ and $u = \nabla \mathcal{J}_L(\rho)$. It can be then deduced
from \eqref{eq:defG} that
\begin{equation}
\forall \mathbf{r}
\in \{1,\ldots,Y\}\times \{1,\ldots,X\},\qquad 
\vect{u}(\mathbf{r}) = 2\vect{S}^{\hermit}(\mathbf{r})\vect{\Psi}^{-1}\left( \vect{S}(\mathbf{r})
\boldsymbol{\rho}(\mathbf{r})-\vect{d}(\mathbf{r})\right).\nonumber
\end{equation}
where the vector $\vect{u}(\mathbf{r})$ is defined from $u$ in the same way as $\overline{\boldsymbol{\rho}}(\mathbf{r})$ is defined from $\overline{\rho}$ in \eqref{eq:defvrho}. 
Furthermore, for every $\zeta' \in \mathbb{C}^K$,
\begin{equation}
\|\nabla \mathcal{J}_{LT}(\zeta)- \nabla \mathcal{J}_{LT}(\zeta')\|
\le \|T\| \|u-u'\|
\label{eq:lips1}
\end{equation}
where  $u' = \nabla \mathcal{J}_L(\rho')$ and $\rho' = T^*\zeta'$. We have then
\begin{align}
\|u-u' \|^2 &= \sum_{\mathbf{r}
\in \{1,\ldots,Y\}\times \{1,\ldots,X\}} 
\|\vect{u}(\mathbf{r})-\vect{u'}(\mathbf{r})\|^2\nonumber\\
& =4\sum_{\mathbf{r}
\in \{1,\ldots,Y\}\times \{1,\ldots,X\}} 
\|\vect{S}^{\hermit}(\mathbf{r})\vect{\Psi}^{-1}\vect{S}(\mathbf{r})\left(
\boldsymbol{\rho}(\mathbf{r})-\vect{\rho'}(\mathbf{r})\right)\|^2\nonumber\\
& \le 4\sum_{\mathbf{r}
 \in \{1,\ldots,Y\}\times \{1,\ldots,X\}}  \theta_{\mathbf{r}}^2
\|\boldsymbol{\rho}(\mathbf{r})-\boldsymbol{\rho}'(\mathbf{r})\|^2\nonumber\\
& \le 4\theta^2 \|\rho-\rho'\|^2\nonumber\\
& \le 4\theta^2 \|T\|^2 \|\zeta-\zeta'\|^2.
\label{eq:lips2}
\end{align}
Altogether, \eqref{eq:lips1} and \eqref{eq:lips2} yield
\begin{equation}
\|\nabla \mathcal{J}_{LT}(\zeta)- \nabla \mathcal{J}_{LT}(\zeta')\|
\le 2\theta\|T\|^2 \|\zeta-\zeta'\|
\end{equation}
which shows that $\mathcal{J}_{LT}$ has a Lipschitz continuous gradient with
constant $2\theta \|T\|^2$.

Based on these observations and the fact that, 
\begin{multline}
\forall \zeta=\big((\zeta_{a,k})_{1 \le k\le K_{j_\mathrm{max}}}, (\zeta_{o,j,k})_{1 \le j \le j_\mathrm{max}, 1\le k \le K_j}\big),\\
\prox_{\gamma _{n} \mathcal{J}_P} \zeta =
\big((\prox_{\gamma_n \Phi_a}\zeta_{a,k})_{1 \le k\le K_{j_\mathrm{max}}}, (\prox_{\gamma_n \Phi_{o,j}}\zeta_{o,j,k})_{1 \le j \le j_\mathrm{max}, 1\le k \le K_j}\big),
\label{eq:proxJ2}
\end{multline} 
the sequence $(\zeta^{(n)})_{n>0}$ built by Algorithm \ref{algo:wt}
can be rewritten under the more classical forward-backward iterative form \cite{chaux_07}:
\begin{equation}
\zeta^{(n+1)} = \zeta^{(n)}+\lambda_{n}\Big(
\mathrm{prox}_{\gamma _{n} \mathcal{J}_P}
\big(\zeta^{(n)}-\gamma_{n}\nabla \mathcal{J}_{LT}(\zeta^{(n)})\big)-\zeta^{(n)}\Big)
\label{eq:FW}
\end{equation}
and, due to the Lipschitz differentiability of $\mathcal{J}_{LT}$, the convergence
of the algorithm is secured under Assumption~\ref{as:step} (see \cite{Chaux_C_07,Combettes_05}). Furthermore, since $\mathcal{J}_P$ is strongly convex
with modulus $\vartheta_1$, we have (see \cite{pustelnik_08} and references therein)
\begin{equation}
\forall n >0,\qquad \|\zeta^{(n)} -\widehat{\zeta}\| \le 
\Big(1-\frac{\underline{\lambda}\underline{\gamma}\vartheta_1}{1+\underline{\gamma}\vartheta_1}\Big)^{n-1} \|\zeta^{(1)} -\widehat{\zeta}\|
\label{eq:linconv}
\end{equation}
where $\underline{\gamma}=\inf_{n>0} \gamma_n$ and 
$\underline{\lambda} = \inf_{n>0} \lambda_n$. This proves that $(\zeta^{(n)})_{n>0}$ converges linearly to $\widehat{\zeta}$.

\footnotesize
\bibliographystyle{elsarticle-num}
\bibliography{biblio_chaari}

\end{document}